\documentclass[10pt,a4paper]{article}
\usepackage{a4wide}
\usepackage{amsbsy,amscd,amsfonts,amssymb,amstext,amsmath,latexsym,theorem}

\usepackage{ifthen}

\newboolean{isArxivVersion}
\setboolean{isArxivVersion}{true}

\usepackage{listings}
\lstdefinelanguage{gap} 
        {morekeywords={and,break,continue,do,elif,else,end,false,fi,for,
                function,if,in,local,not,od,or,repeat,return,then,true,until,while}, 
        sensitive, 
        morecomment=[l]\#, 
        morestring=[b]", 
        morestring=[b]', 
} 
\usepackage[all]{xy} 
\mathchardef\tnode="020E 
 
\def\arc{
  \hbox{\kern -0.15em 
  \vbox{\hrule width 4.5em height 0.6ex depth -0.5 ex} 
  \kern -0.33em}} 
 
\def\darc{
  \rlap{\lower0.2ex\arc}{\raise0.2ex\arc}} 
 
\def\stroke#1{
  \kern 0.05em 
  \rlap\arc{{\textstyle{#1}}\atop\phantom\arc} 
  \kern -0.22em} 
 
\def\dstroke#1{
  \kern 0.05em 
  \rlap\darc{{\tiny\textstyle{#1}}\atop\phantom\darc} 
  \kern -0.22em} 
 
\def\centerscript#1{
  \setbox0=\hbox{$\tnode$} 
  \hbox to \wd0{\hss$\scriptstyle{#1}$\hss}} 
 
\def\node{
  \def\super{} 
  \def\sub{} 
  \futurelet\next\dolabellednode} 
 
  \let\sp=^ 
  \let\sb=_ 
 
  \def\dolabellednode{%
    \ifx\next\sb\let\next\getsub 
    \else 
      \ifx\next\sp\let\next\getsuper 
      \else\let\next\donode 
      \fi 
    \fi 
    \next} 
 
  \def\getsub_#1{\def\sub{#1}\futurelet\next\dolabellednode} 
  \def\getsuper^#1{\def\super{#1}\futurelet\next\dolabellednode} 
 
  \def\donode{%
   \rlap{$\mathop{\phantom\tnode}\limits_{\centerscript{\sub}}^{\centerscript{\super}}$}\tnode}

\def\varcdn{
  \kern -0.03em\vbox{\kern -0.5ex 
  \hbox to \wd0{\hss\vrule width 0.04em depth 5.8ex\hss} 
  \kern -0.3ex  \hbox{$\tnode$}}}

\newcommand{\Theorem}{Theorem} 
\newcommand{\Curtistits}{Curtis-Tits Theorem} 
\newcommand{\Titslem}{Tits' Lemma} 
\newcommand{\Proposition}{Proposition} 
\newcommand{\Lemma}{Lemma} 
\newcommand{\Corollary}{Corollary} 
\newcommand{\Program}{Program} 
\newcommand{\Definition}{Definition} 
\newcommand{\Remark}{Remark} 
\newcommand{\Example}{Example} 
\newcommand{\Examples}{Examples} 
\newcommand{\Consequence}{Consequence} 
\newcommand{\Fact}{Fact} 
\newcommand{\Notation}{Notation} 
\newcommand{\Problem}{Problem} 
\newcommand{\Conjecture}{Conjecture} 
\newcommand{\Technique}{{\sc Technique}} 
\theoremstyle{break} 
\newtheorem{theorem}{\Theorem}[section]

\newtheorem{proposition}[theorem]{\Proposition} 
\newtheorem{lemma}[theorem]{\Lemma} 
\newtheorem{corollary}[theorem]{\Corollary} 
{\theorembodyfont{\sc} 
}

\theoremstyle{plain} 
{\theorembodyfont{\rmfamily} 
 
\newtheorem{definition}[theorem]{\Definition}

\newtheorem{remark}[theorem]{\Remark} 
\newtheorem{example}[theorem]{\Example} 
 
\theoremstyle{break}

} 
{\theorembodyfont{\ttfamily} 
\theoremstyle{break} 
 
} 
\newenvironment{proof}{\noindent {\sl Proof. }}{\hfill \pend \br}

\newcommand{\mc}[1]{\mathcal{#1}}

\newcommand{\pend}{\hfill $\Box$} 
\newcommand{\nend}{\hfill $\blacksquare$} 
 
\newcommand{\br}{\vspace{10pt}} 
 
\newcommand{\lb}{\left\{} 
 
\newcommand{\<}{\left<} 
\newcommand{\rb}{\right\}} 
 
\renewcommand{\>}{\right>} 
\newcommand{\gen}[1]{\left< #1 \right>}

\newcommand{\typ}{{\rm typ}} 
 
\newcommand{\Aut}{{\rm Aut}} 
 
\newcommand{\G}{{\mc{G}}}

\newcommand{\PG}{{\rm PG}} 
 
\newcommand{\Sym}{\rm Sym}

\renewcommand{\Sym} {{\rm Sym}}

\newcommand\A{{\cal A}}

\newcommand\K{\mathbb{K}}

\renewcommand\hat{\widehat}

\newcommand{\sfG}{\mathsf{G}} 
\newcommand{\sfQ}{\mathsf{Q}} 
\newcommand{\sfH}{\mathsf{H}} 
\newcommand{\sfW}{\mathsf{W}} 

\newcommand{\abs}[1]{\lvert#1\rvert}        
 
\begin{document} 
\title{{\bf Intransitive geometries and fused amalgams}
\ifthenelse{\boolean{isArxivVersion}}{\\Extended arXiv version}{}
}
\author{Ralf Gramlich \and Max Horn \and Antonio Pasini \and Hendrik Van Maldeghem}
\date{\today}
\maketitle

\begin{abstract}
\noindent We study geometries that arise from the natural $\sfG_2(\mathbb{K})$ action on the geometry of one-dimensional subspaces, of nonsingular two-dimensional subspaces, and of nonsingular three-dimensional subspaces of the building geometry of type $C_3(\mathbb{K})$ where $\mathbb{K}$ is a perfect field of characteristic $2$. One of these geometries is intransitive in such a way that the non-standard geometric covering theory from \cite{Gramlich/Maldeghem} is not applicable. 
In this paper we introduce the concept of fused amalgams in order to extend the geometric covering theory so that it applies to that geometry. This yields an interesting new amalgamation result for the group $\sfG_2(\mathbb{K})$.
\end{abstract} 
 
\section{Introduction} 
 
Tits' lemma \cite{Tits:1986} (see also \cite[Lemma 5]{Pasini:1985} or \cite[Theorem 12.28]{Pasini:1994}) provides a geometric way to prove that certain groups can be identified as universal enveloping groups of certain amalgams. More
precisely, the universal enveloping group of the amalgam of parabolic subgroups of a group $G$ acting flag-transitively on a geometry
$\Gamma$ equals $G$ if and only if
$\Gamma$ is simply connected. Obviously, this technique to compute amalgams of groups is limited. In order
to include more amalgams with this geometric technique, one can generalise Tits' lemma to intransitive
geometries. Roughly speaking, two difficulties have to be overcome in this generalisation process: (1) the
reconstruction of the geometry from the various parabolic subgroups, (2) finding the right amalgam and getting
control over its universal enveloping group by using the simple connectivity of the geometry. 

In \cite{Gramlich/Maldeghem} a theory was established for geometries with an automorphism group admitting possibly more than one vertex orbit per type.
For the
reconstruction of the geometry from the stabiliser data the authors used a result of Stroppel \cite{Stroppel:1993}. 
Their work was motivated by the
geometries arising from non-isotropic elements with respect to an orthogonal polarity in projective space. Another sporadic amalgam, considered by Hoffman and Shpectorov
\cite{Hoffman/Shpectorov}, related to the group $\mathsf{G}_2(3)$ could be handled very elegantly with that
theory. Yet another application of the new covering theory is a local characterisation of the group $\mathsf{SL}_{n+1}(\mathbb{F}_q)$ via centralisers of root subgroups using the local characterisation of the graph on incident point-hyperplane pairs obtained in \cite{Gramlich/Pralle}.

In general, it seems that non-standard amalgams related to exceptional groups of Lie type cannot be treated
with Tits' original lemma \cite{Tits:1986}. But also the intransitive theory developed in
\cite{Gramlich/Maldeghem} often falls short, as it requires that 
\begin{description}
\item[$(\ddagger)$] for every flag $F$ of rank two or 
three of $\Gamma$, the action of $G$ on the orbit $G.F$ is flag-transitive. 
\end{description}

In the present paper, we consider some rather natural amalgams
related to Dickson's groups of type $\mathsf{G}_2$. For some of them, the existing theory suffices to get control
over the universal enveloping group. For others, we need to modify the theory. This will lead us to \emph{fused
amalgams}. Roughly speaking, fused amalgams occur when the corresponding group acts intransitively on the set
of maximal flags of the corresponding geometry as in \cite{Gramlich/Maldeghem} and Stroppel's reconstruction of
the geometry fails because Property~$(\ddagger)$ above is not satisfied. In our example, the group acts
transitively on each type of vertex (and there are three types), but there are two orbits on the set of
chambers. As a result, there seems to be no purely group-theoretic way to reconstruct the geometry. Instead, we use
the properties of the diagram to settle incidences that cannot be recovered by the group. In particular we exploit the
fact that the diagram is a string of length three, and that hence the residues of the elements belonging to the
middle node are generalised digons.
We only consider geometries belonging to tree diagrams, because circuits
introduce ambiguity on how to define incidence in residues that are generalised digons. The main covering theoretic results of this article are the Reconstruction Theorem \ref{isomorphism theorem2} and the Covering Theorem \ref{mainthmcover}. The main applications of this covering theory contained in this paper are the simple connectivity result Theorem \ref{simplemain} and the amalgamation result Theorem \ref{amalg}  

The paper is organised as follows. In Section~\ref{s2} we define the geometries that are relevant for the rest
of the paper and in particular for the amalgams that we will consider in Section~\ref{s3}. In Section~\ref{s4}, we develop a theory of intransitive geometries and amalgams, which
we call \emph{fused amalgams}, that allows us to tackle the amalgam described in Section~\ref{s3}. We immediately apply the theory to our
situation. In the final Section~\ref{s5}, we establish the simple connectivity for the investigated geometries.
 
\section{Some geometries related to $\mathsf{G}_2$}\label{s2}

\subsection{The split Cayley hexagon and its properties} \label{2.1}

We consider the Chevalley group $\sfG_2(\K)$, with $\K$ any (commutative) field. Naturally associated with each
Chevalley group is a \emph{building}, in the sense of Tits \cite{Tits:1974}. In the case of $\sfG_2(\K)$, this
building is a bipartite graph, which is the incidence graph of a pair of generalised hexagons (and we may
freely consider one of those by choosing one of the bipartition classes as set of points, and the other class
as set of lines). We use the standard notions from incidence geometry and the theory of building geometries, like distances between elements and opposition of elements, cf.\ \cite{Pasini:1994} and \cite{Tits:1974}. We recall that a \emph{generalised hexagon} is a point-line incidence geometry with the
properties that
 
\begin{itemize}
\item[(GH1)] every two elements (points or lines) are contained in an ordinary hexagon (i.e., a cycle of 12
distinct consecutively incident elements), and

\item[(GH2)] there are no ordinary $n$-gons for $n<6$.
\end{itemize}

The generalised hexagon related to $\sfG_2(\K)$ is called the \emph{split Cayley hexagon} and can be
represented on the parabolic quadric $\sfQ(6,\K)$, which is a nondegenerate quadric in $\PG(6,\K)$ of (maximal)
Witt index $3$. The points of the hexagon are all points of $\mathsf{Q}(6,\K)$, while the lines (the
\emph{hexagon lines}) are only some well-chosen lines on $\mathsf{Q}(6,\K)$. If $\mathsf{Q}(6,\K)$ has the
standard equation $X_0X_4+X_1X_5+X_2X_6=X_3^2$, then a line on $\mathsf{Q}(6,\K)$ with Grassmannian coordinates
$(p_{01},p_{02},\ldots,p_{06},p_{12},p_{13},\ldots,p_{56})$ is a hexagon line if and only if
$p_{12}=p_{34}$, $p_{56}=p_{03}$, $p_{45}=p_{23}$, $p_{01}=p_{36}$, $p_{02}=-p_{35}$ and $p_{46}=-p_{13}$.  This line set
has the following properties (see e.g.~\cite{Maldeghem:1998}).

\begin{enumerate}
\item The set of lines of $\sfH(\K)$ through a fixed point fills up a projective plane on $\mathsf{Q}(6,\K)$.
We call such a plane a \emph{hexagonal plane}. 
\item  Any plane of $\mathsf{Q}(6,\K)$ that contains
at least one hexagon line is a hexagonal one. Any other plane of $\mathsf{Q}(6,\K)$ will be called
an \emph{ideal plane}. 
\item Every line of $\sfQ(6,\K)$ that does not belong to the hexagon is contained in a
unique hexagonal plane. We call such a line an \emph{ideal line}. The point of the corresponding hexagonal
plane that is the intersection of all hexagon lines in that plane is called the \emph{ideal centre} of the
ideal line. 
\item The ideal centres of all ideal lines of an ideal plane $\pi$ form again an ideal plane
$\pi'$, which, together with $\pi$, generates a hyperplane $H$ of $\PG(6,\K)$ that intersects $\sfQ(6,\K)$ is a
non-degenerate (hyperbolic) quadric of Witt index $3$. The point set of $\pi\cup\pi'$ is the point set of a
non-thick ideal subhexagon of $\sfH(\K)$. The hyperplane $H$ is called a \emph{hyperbolic hyperplane}. Every hyperplane that intersects $\sfQ(6,\K)$ in a hyperbolic quadric arises in this way. In particular, every hyperplane $H$ that contains a plane of $\sfQ(6,\K)$ is either a tangent hyperplane or a hyperbolic
hyperplane. The former does not contain disjoint planes of $\sfQ(6,\K)$.
\end{enumerate}
 
If $\K$ has characteristic two, and $\K$ is perfect (which means that the mapping $x\mapsto x^2$ is
surjective), then the projection of the point set of $\sfQ(6,\K)$ from the point $(0,0,0,1,0,0,0)$ onto the
hyperplane $\PG(5,\K)$ with equation $X_3=0$ embeds $\sfQ(6,\K)$ bijectively onto a symplectic space
$\mathsf{W}(5,\K)$, so that we also obtain an embedding of $\sfH(\K)$ into $\mathsf{W}(5,\K)$. The lines of $\sfQ(6,\K)$ are projected onto totally isotropic lines with respect to the
corresponding symplectic polarity (we will call such lines \emph{symplectic lines}), cf.\ 2.4.14 of
\cite{Maldeghem:1998}. The lines of $\PG(5,\K)$ that are not symplectic will be called \emph{non-symplectic
lines}. The projection of hexagonal planes and ideal planes will be called \emph{hexagonal} and \emph{ideal},
respectively. Likewise, the projection of hexagon and ideal lines will be called {\em hexagon} and {\em ideal}, respectively; both are symplectic lines. A \emph{nonsingular} plane is a plane of $\PG(5,\K)$ in which the non-symplectic lines form a
dual affine plane. A \emph{special} nonsingular plane is a nonsingular plane containing a hexagon line.
 
The above properties also translate to the situation in $\PG(5,\K)$, when the characteristic of $\K$ is equal
to two. For instance, every ideal line is contained in a unique hexagonal plane and the ideal centre is not
contained in that ideal line.
 
Moreover, we have the following:
\begin{enumerate}
\setcounter{enumi}{4}
\item Let $l$ be a non-symplectic line of $\mathrm{PG}(5,\K)$. Then the set of hexagon lines at hexagon-distance three from all points of $l$ form a distinguished regulus $\mathcal{R}$ of a hyperbolic quadric $\sfQ(3,\K)$ in the orthogonal space $l^\perp$ of $l$ with respect to the symplectic polarity. This follows immediately from the regulus property (see \cite{Ronan:1980}) and the fact that opposition of points in the hexagon corresponds to non-perpendicularity in the symplectic polar space. Moreover, every pair of opposite lines (in the hexagon) is contained in such a regulus.
\end{enumerate}

\subsection{Some geometries} \label{gamma}
 
We now consider four different infinite classes of geometries $\Gamma_0$, $\Gamma_1$, $\Gamma_2$, $\Gamma_3$, all of rank 3 and with type set $\{1,2,3\}$.
We call elements of type 1 \emph{points}, of type 2 \emph{lines}, and of type 3 \emph{planes}.

To define $\Gamma_0$, $\Gamma_1$, $\Gamma_2$ let $\K$ be perfect and of characteristic two. The geometries $\Gamma_i$, $i\in\{0,1,2\}$, have as set of
elements of type 1 the set of points of $\PG(5,\K)$, and as set of elements of type 2 the set of non-symplectic
lines of $\PG(5,\K)$, with natural incidence. The elements of type 3 of the geometries $\Gamma_0$, $\Gamma_1$,
and $\Gamma_2$ are all nonsingular planes, all nonspecial nonsingular planes, and all special nonsingular
planes, respectively. Incidence between elements of type 2 and 3 is natural, and incidence between elements of
type 1 and 3 is given by the following rule: $p$ is incident with $\pi$ if and only if there is a type 2 element $l$
incident with both.
The geometry $\Gamma_0$ is flag-transitive for the symplectic group $\sf{S}_6(\K)$ and has been considered by Cuypers \cite{Cuypers} and Hall \cite{Hall} and, more recently, by 
Blok and Hoffman \cite{Blok/Hoffman}; see also \cite{Gramlich:2004}.

The geometry $\Gamma_0$
is in a certain sense a {\em join} of the geometries $\Gamma_1$ and $\Gamma_2$. Indeed, the point-line truncations of
$\Gamma_0$, $\Gamma_1$, and $\Gamma_2$ coincide, while the plane set of $\Gamma_0$ consists of the disjoint union of the
plane sets of $\Gamma_1$ and $\Gamma_2$. 
 
The following gives a relation between covers of two connected rank $3$ geometries $\Delta_1$, $\Delta_2$ having 
identical point-line truncations and the join $\Delta$ of $\Delta_1$ and $\Delta_2$. For $i = 1, 2$ let $\tilde 
\Delta_i$ be the universal cover of $\Delta_i$ and let $t_i$ be number  of layers of the covering projection 
from $\tilde \Delta_i$ to $\Delta_i$. Furthermore, let $\tilde \Delta$ be the universal cover of $\Delta$ and 
let $t$ be number of layers of the covering projection from   $\tilde \Delta$ to $\Delta$. (We refer the reader to \cite{Seifert/Threlfall:1934} for a thorough introduction to the covering theory of simplicial complexes, in particular for the definition of a universal cover of a connected pure simplicial complex.)
 
\begin{proposition} \label{join}
Assume that the planes of $\Delta_1$ and $\Delta_2$ are connected. If, for some $i \in \{ 1, 2 \}$, we have $t_i < \infty$, then $t | t_i$. 
\end{proposition} 
 
\begin{proof} 
Let $\pi : \tilde \Delta \to \Delta$ be a universal covering of $\Delta$ and, for $i = 1, 2$, let $\overline{\Delta}_i$ be  the pre-image of $\Delta_i$ under $\pi$. Since the planes of $\Delta_i$ are connected and since $\Delta_i$ and $\Delta$ have identical connected point-line truncations, the preimage $\overline{\Delta}_i$ is connected, so $\pi$ induces a covering from $\overline{\Delta}_i$ to $\Delta_i$. Hence $t|t_i$, if $t_i$ is finite.
\end{proof} 
 
\begin{corollary} 
Assume that the planes of $\Delta_1$ and $\Delta_2$ are connected. 
\begin{enumerate} 
\item If one of $\Delta_1$, $\Delta_2$ is simply connected, then $\Delta$ is simply connected. 
\item Suppose $t_1, t_2 < \infty$. Then $t|\gcd(t_1,t_2)$. In particular, $\Delta$ is simply connected, if $t_1$ and $t_2$ are coprime. 
\end{enumerate} 
\end{corollary} 
 
\begin{remark}
The join of two geometries $\Delta_1$ and $\Delta_2$ can be simply connected, even if $\Delta_1$ and 
$\Delta_2$ are isomorphic and admit infinite universal covers. A nice example for this behaviour is the geometry studied in \cite{Hoffman/Shpectorov}, which in fact also occurs in \cite{Aschbacher/Smith}, \cite{Wispelaere/Huizinga/Maldeghem:2005}, \cite{Kantor:1985} in different guise. 
  In 
\cite{Hoffman/Shpectorov} Hoffman and Shpectorov study an amalgam 
of maximal subgroups of $\hat G = \Aut(\mathsf{G_2}(3))$ given by 
a certain choice of subgroups $\hat L =  2^3 \cdot \mathsf{L}_3(2) : 2$, 
$\hat N  =  2^{1+4}_+.(S_3 \times S_3)$, $M  =  \mathsf{G_2}(2) = \mathsf{U}_3(3):2$ 
which corresponds to an amalgam of subgroups of $G = G_2(3)$ given 
by 
$L  =   \hat L \cap G = 2^3 \cdot \mathsf{L}_3(2)$, $N  =  \hat N \cap G = 2^{1+4}_+.(3 \times 3).2$, $M  =  \mathsf{G_2}(2) = \mathsf{U}_3(3):2$, $K  =  eMe^{-1}$ for $e \in O_2(\hat L) \backslash 
O_2(L)$. The groups 
$\hat{G}_1  =  \hat L$, $\hat{G}_2  =  \hat N$, $\hat{G}_3  =  M$ 
define a flag-transitive coset geometry $\Gamma$ of rank three for 
$\hat G = \Aut(\mathsf{G_2}(3))$, which is simply connected by 
\cite{Hoffman/Shpectorov}. The subgroup $G = \mathsf{G_2}(3)$ of 
$\hat G$ does not act flag-transitively on $\Gamma$. Nevertheless, the 
groups 
$G_p  =  L$, 
$G_l  =  N$, 
$G_{\pi_1} =  M$,
$G_{\pi_2}  =  K$  
define an intransitive coset geometry of rank three for $G = 
\mathsf{G_2}(3)$ satisfying Property $(\ddagger)$ from the introduction, which is isomorphic to $\Gamma$ by 
\cite{Hoffman/Shpectorov} and, hence, simply connected, so that non-standard covering theory as in \cite{Gramlich/Maldeghem} is applicable.

The coset geometries $(G_p,G_l,G_{\pi_1},*)$ and $(G_p,G_l,G_{\pi_2},*)$ are isomorphic to the GAB | Geometry that is Almost a Building | studied in \cite[Table 1, Example 4]{Aschbacher/Smith}, in \cite[Section 6.1]{Wispelaere/Huizinga/Maldeghem:2005}, and in \cite{Kantor:1985} with diagram
$$\node\stroke{6}\node\stroke{}\node.$$ 
This GAB is very far from being simply connected. In fact, by \cite{Kantor:1985} the amalgam of $L$, $N$, $M$ admits the group $G_2(\mathbb{Q}_2)$ as a universal enveloping group | while by \cite{Hoffman/Shpectorov} the amalgam of $\hat L$, $\hat N$, $M$ admits the group $\Aut(\mathsf{G_2}(3))$ as its universal enveloping group. We conjecture that Kantor's description \cite{Kantor:1985} of the universal cover of $(G_p,G_l,G_{\pi_1},*) \cong (G_p,G_l,G_{\pi_2},*)$ can be used to give an alternative proof of the simple connectivity of the join of $(G_p,G_l,G_{\pi_1},*)$ and $(G_p,G_l,G_{\pi_2},*)$ by studying those quotients of the group $G_2(\mathbb{Q}_2)$ that admit an involutory outer automorphism.
However, the combinatorial simple connectivity proof given by Hoffman and Shpectorov \cite{Hoffman/Shpectorov} is short and clear and likely to be shorter than any group-theoretic proof of simple connectivity.
\end{remark}

Concerning the fourth class of geometries, let $\K$ be any field. Then the  rank 3 geometry $\Gamma_3$ consists of the points of the split Cayley hexagon $\sfH(\K)$, the ideal lines,
and the ideal planes, with natural incidence. The amalgam and corresponding geometry
$\Gamma_3$ considered here has also been treated by Baumeister, Shpectorov and Stroth in
an unpublished manuscript \cite{Baumeister/Shpectorov/Stroth}. We have found an independent proof which we include here so that a proof of this fact is made available in the literature. Moreover there exists a result \cite{Shpectorov} by Shpectorov dealing with the simple connectivity of hyperplane complements in arbitrary dual polar spaces with line size at least five, thus independently implying simple connectivity of $\Gamma_3$, but only for $|\K| \geq 4$.

\bigskip
In this article we prove the following results:

\begin{theorem} \label{simplemain}
\begin{enumerate}
\item The geometry $\Gamma_0$ is simply connected.
\item The geometry $\Gamma_1$ is flag-transitive. Moreover, it is simply connected, whenever $|\K|>2$.
\item The geometry $\Gamma_2$ is simply connected.
\item The geometry $\Gamma_3$ is simply connected.
\end{enumerate}
\end{theorem}

\begin{proof}
\begin{enumerate}
\item This is proved in Proposition \ref{gamma1}, also \cite{Blok/Hoffman} or Proposition \ref{join} plus Proposition \ref{gamma3}.
\item See Propositions \ref{gamma2t} and \ref{gamma2}.
\item Cf.\ Proposition \ref{gamma3}.
\item This follows from Proposition \ref{prop3}, also \cite{Baumeister/Shpectorov/Stroth}, or \cite{Shpectorov} for $|\K| \geq 4$.
\end{enumerate}
\end{proof}

\subsection{Amalgams for $\Gamma_2$}\label{ssGamma2} \label{s3}

The geometries $\Gamma_0$ and $\Gamma_3$ have been extensively studied. Moreover, the geometry $\Gamma_1$ is flag-transitive, cf.\ Proposition \ref{gamma2t}, so that classical covering theory applies. Hence we concentrate on the amalgam of parabolics given by the $\sf{G}_2(\mathbb{K})$ action on $\Gamma_2$.

Let $p, l, \pi$ be a chamber of $\Gamma_2$ and denote by $G_p$, $G_l$, $G_\pi$, $G_{p,l}$, etc., the respective stabilisers. We now collect information about $\Gamma_2$ and these stabilisers, most of which are based on the following proof of the flag-transitivity of $\Gamma_1$.

\begin{proposition} \label{gamma2t}
The action of $\mathsf{G}_2(\K)$ on the geometry $\Gamma_1$ is flag-transitive. 
\end{proposition} 
 
\begin{proof} Set $G:=\mathsf{G}_2(\K)$. 
 All non-symplectic lines are determined by two opposite points of $\sfH(\K)$. The fact that $G$ 
acts transitively on pairs of opposite points of $\sfH(\K)$ (see e.g.~Chapter 4 of \cite{Maldeghem:1998}) implies that $G$ acts transitively on the point-line pairs of $\Gamma_1$. So we may fix such a point-line pair 
$(x,l)$ and it suffices to prove that $H:=G_{x,l}$ acts transitively on the planes of $\Gamma_1$ containing $l$. The group $H$ 
stabilises the pole $\Sigma$ of $l$ with respect to the symplectic polarity related to the $C_3$ building $\sfW(5,\K)$. The pole $\Sigma$ is a projective $3$-space, so that
every plane $\pi$ of $\PG(5,\K)$ containing $l$ meets $\Sigma$ in a unique point $x_\pi$. Viewed in $\sfH(\K)$, 
the space $\Sigma$ is determined by the lines at distance three from all the points of $l$. These lines form a 
distinguished regulus $\mathcal{R}$ of a hyperbolic quadric $\sfQ(3,\K)$ in $\Sigma$, cf.\ Section \ref{2.1}, item (v). Clearly, $\pi$ is nonsingular. Also, it is easy to see that $\pi$ contains a hexagon line if and only if $x_\pi$ is contained in 
the quadric $\mathsf{Q}(3,\K)$. Let $\mathcal{R}'$ be the complementary regulus of $\mathcal{R}$ on 
$\mathsf{Q}(3,\K)$. Then every point $y$ on $l$ uniquely determines a line $m$ of $\mathcal{R}'$ by the fact 
that all hexagon lines through $y$ meet $m$. Now, the stabiliser in $G$ of $\mathsf{Q}(3,\K)$ contains the 
group $\mathsf{L}_2(\K)\times\mathsf{L}_2(\K)$. Hence the assertion is equivalent with saying that in $\Sigma$, 
the group $\sf{L}_2(\K)\times \sf{L}_2(\K)$ stabilising the hyperbolic quadric $\mathsf{Q}(3,\K)$ acts transitively on the 
pairs $(p,k)$, where $p$ is a point off the quadric, and $k$ is a line of a fixed regulus of the quadric, which 
is a true statement as one can easily verify. 
\end{proof} 

\begin{proposition} \label{properties}
The action of $\sfG_2(\K)$ on the geometry $\Gamma_2$ is transitive on the incident point-line pairs and transitive on the incident line-plane pairs, but it is intransitive on the incident point-plane pairs and has two incident point-plane orbits instead. Moreover, the stabiliser of an incident line-plane pair $(l,\pi)$ has two orbits on the points incident to $l$.
\end{proposition}

\begin{proof}
\begin{enumerate}
\item {Point-line-transitivity:} The point-line truncations of $\Gamma_1$ and $\Gamma_2$ coincide (see the discussion before Proposition \ref{join}) and $\sfG_2(\mathbb{K})$ is flag-transitive on $\Gamma_1$ by Proposition \ref{gamma2t}, so that $\sfG_2(\mathbb{K})$ acts transitively on the incident point-line pairs of $\Gamma_2$.
\item {Point-plane-intransitivity:} The planes of $\Gamma_2$ are those rank 2 planes of the $C_3$ building geometry $\sfW(5,\mathbb{K})$ which contain a (unique) hexagon line. Therefore the plane stabiliser $G_\pi$ has to fix this hexagon line and consequently cannot map a point on that line onto a point in the plane not on the line. Hence $G_\pi$ is not transitive on the set of points incident with $\pi$, whence $\Gamma_2$ is not point-plane-transitive.
\item {Line-plane-transitivity:} Let $l$ be a line of $\Gamma_2$ as in the proof of Proposition \ref{gamma2t}. A plane $\pi$ of $\mathrm{PG}(5,\mathbb{K})$ containing $l$ has rank two with respect to the symplectic form, and intersects the pole $\Sigma$ of $l$ with respect to the symplectic form in a point $x_\pi$. As in the proof of Proposition \ref{gamma2t}, the space $\Sigma$ carries the structure of a $\mathsf{Q}(3,\mathbb{K})$, and $\pi$ contains a hexagon line if and only if $x_\pi$ is contained in $\mathsf{Q}(3,\mathbb{K})$. Line-plane-transitivity now is a consequence of line-transitivity and transitivity of $G_l$ on the points of  $\mathsf{Q}(3,\mathbb{K})$.
\item Two orbits: Denote the hexagon line contained in $\pi$ by $h$. The pole $l^\perp$ of $l$ with respect to the symplectic polarity contains a regulus of hexagon lines, cf.\ Section \ref{2.1}, item (v). The map sending a point of $l$ onto the set of points of $\mathcal{R}$ at distance three in $\sfH(\K)$ is a bijection of the points of $l$ onto the lines of the complementary regulus of $\mathcal{R}$. Since the pointwise stabiliser $H$ in $\sfG_2(\K)$ of $\mathcal{R}$ acts two-transitively on the complementary regulus (this follows from the Moufang property and the regulus condition, cf.\ \cite{Ronan:1980}, \cite[Proposition 4.5.11]{Maldeghem:1998}) and since $H$ stabilises the line $l$, we see that $H_x$ acts transitively on $l \backslash \{ x \}$ where $x = l \cap h$. Hence $\sfG_2(\K)$ has two orbits on the incident point-plane pairs as well.
\end{enumerate}
\end{proof}

%

\section{Fused amalgams and intransitive geometries}\label{s4}

The nonstandard notions that we will need below were introduced in \cite{Gramlich/Maldeghem}, to which we refer
for more details and results.

\subsection{Diagram coset pregeometries}

\begin{definition}[Diagram Coset Pregeometry] \label{cospre}
Let $I$ be a finite set, let $\Delta = (I,\sim)$ be a tree, and let $(T_i)_{i \in I}$ be a family of pairwise disjoint sets. Also, let $G$ be a group and let
$(G^{t,i})_{t \in T_i, i \in I}$ be a family of subgroups of $G$. Then the {\em diagram coset pregeometry of $G$} with respect to $(G^{t,i})_{t \in T_i, i
\in I}$ equals the pregeometry
$$\left( \left\{(C,t):t\in T_i\mbox{ for some }i\in I, C\in G/G^{t,i}
\right\},*,\typ \right)$$ over $I$ with $\typ (C,t) = i$ if $t\in T_i$, and
\begin{description}
\item[(DCos)] $gG^{t,i} * hG^{s,j}$ if
\begin{itemize}
\item $i = j$ and $t=s$ and $gG^{t,i} \cap hG^{s,j} \neq \emptyset$, \item $i$, $j$ adjacent in $\Delta$ and
$gG^{t,i} \cap hG^{s,j} \neq \emptyset$, or \item $i$, $j$ not adjacent in $\Delta$ and there exists a geodesic
$i = x_0, \ldots, x_k = j$ in $\Delta$ and cosets $g_{t_{x_l},x_l}G^{t_{x_l},x_l}$ with $gG^{t,i} =
g_{t_{x_0},x_0}G^{t_{x_0},x_0}$, $hG^{s,j} = g_{t_{x_k},x_k}G^{t_{x_k},x_k}$ and
$g_{t_{x_l},x_l}G^{t_{x_l},x_l}*g_{t_{x_{l+1}},x_{l+1}}G^{t_{x_{l+1}},x_{l+1}}$.
\end{itemize}
\end{description}
Since the type function is completely determined by the indices, we also denote the coset pregeometry
of $G$ with respect to $(G^{t,i})_{t \in T_i, i \in I}$ by
$$((G/G^{t,i}\times\{t\})_{t \in T_i, i \in I}, *).$$
If the diagram coset pregeometry happens to be a geometry, then $\Delta$ is its basic diagram if and only if at
least one of the residues corresponding to adjacent $i$, $j$ is not a generalised digon.
\end{definition}

\begin{theorem}[inspired by Buekenhout \& Cohen \cite{Buekenhout/Cohen}]
\label{char conn} Let $\abs{I}>1$. The diagram coset geometry $((G/G^{t,i}\times \{t\})_{t \in T_i, i \in I},
*)$ is connected if and only if
$$G= \langle G^{t,i} \mid  i \in I, t \in T_i \rangle.$$
\end{theorem}

\begin{proof}
Suppose that $\Gamma$ is connected. Take $i\in I$ and $t \in T_i$. If $a\in G$, then there is a path
$$1G^{t,i}, a_0 G^{t_0, i_{0}},
a_1 G^{t_1,i_{1}}, a_{2} G^{t_2,i_{2}}, \ldots, a_{m} G^{t_m,i_{m} }, a G^{t,i}$$ in the geometry connecting
the elements $1G^{t,i}$ and $a G^{t,i}$ of $\Gamma$. Extending that path, if necessary, we can assume that the
types $i_j$ and $i_{j+1}$ are adjacent in $\Delta$ for all $j$. Therefore
$$a_k G^{t_k,i_k} \cap a_{k+1} G^{t_{k+1},i_{k+1}} \neq \emptyset,$$
so $$a_k^{-1} a_{k+1}\in G^{t_k,i_k} G^{t_{k+1},i_{k+1}}$$ for $k=0 ,\ldots, m-1$. Hence
$$a = (1^{-1}a_0)( a_0^{-1} a_1 ) \cdots
(a_{m-1}^{-1} a_{m})(a_m^{-1}a) \in  G^{t,i}G^{t_0,i_0} \cdots G^{t_{m-1},i_{m-1}} G^{t_m, i_{m}}G^{t,i},$$ and
so $a\in \langle G^{t,i} \mid  i \in I, j \in T_i \rangle$. The converse is obtained by reversing the above
argument. The only difficulties that can occur are the occasions in which $g_1G^{t_1,i_1} \cap g_2G^{t_2,i_2} \neq
\emptyset$, where $i_1 = i_2$ or $i_1$, $i_2$ not neighbors in $\Delta$. However, this can be remedied by
including some suitable chain of cosets  between $g_1G^{t_1,i_1}$ and $g_2G^{t_2,i_2}$ into the chain of
incidences.
\end{proof}
 
\begin{definition}[Sketch] \label{sketch} 
Let $\Gamma = (X,*,\typ)$ be a geometry over a finite set $I$ whose basic diagram $\Delta$ is a tree, let $G$ 
be a group of automorphisms of $\Gamma$, and let $W \subset X$ be a set of $G$-orbit representatives of $X$. We 
write $$W=\bigcup_{i\in I}W_i$$ with $W_i\subseteq\typ^{-1}(i)$. The {\em sketch of $\Gamma$ with respect to $(G, W, \Delta)$} is the diagram coset geometry 
$$((G/G_w\times\{w\})_{w\in W_i, i\in I},*).$$ 
\end{definition} 
 
Let $\phi : G \to \Sym\ X$ be a group action. Then we denote by $_GX$ the corresponding permutation group, called a {\em $G$-set}.
Two $G$-sets $_GX$ and $_GX'$ are said to be 
{\em equivalent} if there is a bijection $\psi : X \to X'$ such that $\psi \circ \phi(g) \circ \psi^{-1} = 
\phi'(g)$ for each $g\in G$ or, equivalently, $\psi \circ \phi(g) = \phi'(g) \circ \psi$ for all $g\in G$. In 
this case, we shall also say that $_GX$ and $_GX'$ are {\em isomorphic $G$-sets}. 
 
\medskip
Recall also from \cite[Definition 2.2]{Gramlich/Maldeghem} that a {\em lounge} of a geometry $\Gamma = (X,*,\typ)$ over $I$ is a set $W \subseteq X$ of elements such that each subset $V \subseteq W$ for which $\typ_{|V} : V \to I$ is an injection, is a flag. A {\em hall} is a lounge $W$ with $\typ(W) = I$.

\medskip
Finally, recall that, for geometries $\Gamma_1 = (X_1,*_1,\typ_1)$ over $I$ and $\Gamma_2 = (X_2,*_2,\typ_2)$ over $I'$, the {\em direct sum} $\Gamma_1 \oplus \Gamma_2$ is the geometry $(X_1 \sqcup X_2, *_\oplus, \typ_\oplus)$ over $I \sqcup I'$ with ${*_\oplus}_{|X_1 \times X_1} = *_1$ and ${*_\oplus}_{|X_2 \times X_2} = *_2$ and ${*_\oplus}_{|X_1 \times X_2} = X_1 \times X_2$ and ${\typ_\oplus}_{|X_1} = \typ_1$ and ${\typ_\oplus}_{|X_2} = \typ_2$. A geometry $\Gamma$ is said to have the {\em direct sum property}, if for each flag $F$ of $\Gamma$, the residue of $\Gamma$ in $F$ is isomorphic to the direct sum of its truncations to the connected components of its diagram, where the direct sum of more than two geometries is defined iteratively. Note that residual connectivity is a sufficient condition for the direct sum property, see \cite[Theorem 4.2]{Pasini:1994}.

\begin{theorem}[Reconstruction theorem] \label{isomorphism theorem2} 
Let $\Gamma = (X,*,\typ)$ be a geometry over a finite set $I$ with the direct sum property whose basic diagram 
$\Delta$ is a tree. Let $G$ be a group of automorphisms of $\Gamma$. For each $i \in I$ let $$w^i_1, \ldots, 
w^i_{t_i}$$ be $G$-orbit representatives of the elements of type $i$ of $\Gamma$ such that 
\begin{enumerate} 
\item $W := \bigcup_{i \in I} \lb w_1^i, \ldots, w_{t_i}^i \rb$ is a hall and, \item if $V \subseteq W$ is a 
flag, the action of $G$ on the pregeometry over $\typ(V)$ consisting of all elements of the $G$-orbits $G.x$, 
$x \in V$, is transitive on the flags of type $\{ i, j \}$ for all $i, j \in \typ(V)$ corresponding to adjacent 
nodes of the diagram $\Delta$. 
\end{enumerate} 
Then the bijection $\Phi$ between the sketch of $\Gamma$ with respect to $(G, W, \Delta)$ and the pregeometry 
$\Gamma$ given by $$gG_{w_{k}^i} \mapsto gw_k^i$$ is an isomorphism between geometries and an isomorphism between 
$G$-sets. 
\end{theorem} 
 
\begin{proof} 
The isomorphism as $G$-sets is clear from the fundamental theorem of permutation representations as $W$ is a 
transversal with respect to the action of $G$. Therefore let us turn to the isomorphism as geometries. For 
adjacent $i$, $j$ we have $gG_{w_{k_i}^i} \cap hG_{w_{k_j}^j} \neq \emptyset$ if and only if $gw_{k_i}^i * 
hw_{k_j}^j$ by the isomorphism theorem for incidence-transitive geometries. If $i$ and $j$ are non-adjacent, 
then each pair of incident $gw_{k_i}^i * hw_{k_j}^j$ is contained in a chamber of $\Gamma$, hence the basic diagram 
$\Delta$ implies incidence of $gG_{w_{k_i}^i}$, $hG_{w_{k_j}^j}$ in the sketch. If $gw_{k_i}^i$, $hw_{k_j}^j$ 
are not incident, then $gG_{w_{k_i}^i}$, $hG_{w_{k_j}^j}$ cannot be incident by the direct sum property. 
\end{proof} 
 
The direct sum property in the hypothesis of Theorem \ref{isomorphism theorem2} is necessary. Indeed, let 
$\Gamma_1$ and $\Gamma_2$ be isomorphic geometries of rank $4$ with a string basic diagram over the 
type set $\{ 0, 1, 2, 3 \}$. Let $f : \Gamma_1 \to \Gamma_2$ be an isomorphism and glue $\Gamma_1$ to 
$\Gamma_2$ via $f$ restricted to the set of elements of type $1$. Denote the resulting geometry by $\Gamma$. 
Two elements of $\Gamma$ of type distinct from $1$ are incident if and only if they are contributed by the same 
$\Gamma_i$ and are incident in $\Gamma_i$. For $x \in \Gamma_1$ of type $1$ and $y \in \Gamma_i$, we have $x * 
y$ in $\Gamma$ if and only if $x * y$ in $\Gamma_1$ or $f(x) * y$ in $\Gamma_2$. The basic diagram of $\Gamma$ 
also is a string, but the residue of an element of type $1$ does not split into the direct sum of two 
geometries, so $\Gamma$ does not satisfy the direct sum property. Moreover, if the $\Gamma_i$ are 
flag-transitive, then $\Gamma$ is flag-transitive. 
 Altogether, all hypotheses of Theorem \ref{isomorphism theorem2} are satisfied. Nevertheless, $\Gamma$ cannot 
be recovered from its sketch (considered as a diagram coset geometry), because, given an element $x$ of type $1$, Definition \ref{cospre} forces all 
elements of type $0$ incident with $x$ to be incident to all elements of type $2$ incident with $x$, which is 
not the case in $\Gamma$. 
 Of course, $\Gamma$ can be reconstructed in the classical way from its sketch as a flag-transitive geometry, emphasising that 
our reconstruction approach in the present paper is not a generalisation of the classical reconstruction or the 
reconstruction by Stroppel \cite{Stroppel:1993}, cf.\ also \cite{Gramlich/Maldeghem}. 
 
\subsection{Fused amalgams} 
 
In the present paper we will work with the following definition of an amalgam.

\begin{definition}[Amalgam]
Let ${\cal J} = (J,\leq)$ be a finite graded poset with grading function $\tau:J\rightarrow I = \{1,2,...,n\}$ such that every maximal chain has length $n-1$ (namely, it contains an element of every grade). 
 Then an {\em amalgam of shape $\cal J$} is a pair ${\cal A} = ((G_j)_{j\in J}, (\phi_{i,j})_{i<j})$ such that $G_j$ is a group for every $j\in J$ and, for any $i, j \in J$ with $i < j$, the map $\phi_{i,j}:G_i\rightarrow G_j$ is a monomorphism satisfying $\phi_{j,k}\phi_{i,j} = \phi_{i,k}$ for any choice of $i, j, k\in J$ with $i < j < k$.  

The fibers $\tau^{-1}(i)$ are denoted by $J_i$, for $i \in I$.
\end{definition} 
 
\begin{example} \label{2.18} 
In the following diagram we depict an amalgam with $I = \lb 0, 1, 2 \rb$, $J_{0} = \lb 1, 2 \rb$, $J_{1} = \lb 1, 2, 3, 4, 5 \rb$, $J_{2} = \lb 1, 
2, 3, 4 \rb$. The maps $\phi_{i,j}$ are given by arrows and compositions of arrows. 
$$ 
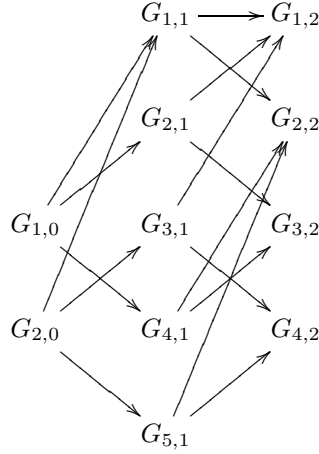
\xymatrix{ 
& G_{1,1} \ar[r] \ar[dr] & G_{1,2} \\ 
& G_{2,1} \ar[ur] \ar[dr]  & G_{2,2} \\ 
G_{1,0} \ar[uur] \ar[ur] \ar[dr] & G_{3,1} \ar[uur] \ar[dr] & G_{3,2} \\ 
G_{2,0} \ar[uuur] \ar[ur] \ar[dr] & G_{4,1} \ar[uur] \ar[ur] & G_{4,2} \\ 
& G_{5,1} \ar[ur] \ar[uuur] & } 
$$ 
 
In terms of stabilisers of a group $G$ acting on a geometry $\Gamma$ with orbit representatives $p$, $l$, 
$\pi_1$, $\pi_2$ this example might concretely arise as 
$$ 
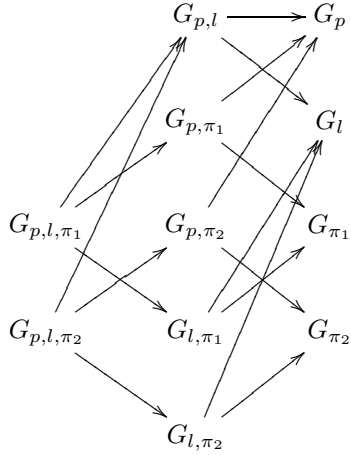
\xymatrix{ 
& G_{p,l} \ar[r] \ar[dr] & G_{p} \\ 
& G_{p,\pi_1} \ar[ur] \ar[dr]  & G_{l} \\ 
G_{p,l,\pi_1} \ar[uur] \ar[ur] \ar[dr] & G_{p,\pi_2} \ar[uur] \ar[dr] & G_{\pi_1} \\ 
G_{p,l,\pi_2} \ar[uuur] \ar[ur] \ar[dr] & G_{l,\pi_1} \ar[uur] \ar[ur] & G_{\pi_2} \\ 
& G_{l,\pi_2} \ar[ur] \ar[uuur] & } 
$$ 
\end{example} 
 
If in the above example $\pi_1$ and $\pi_2$ happen to be contained in the same $G$-orbit, then of course we 
have $G_{\pi_2} = gG_{\pi_1}g^{-1}$ for some $g \in G$. But it may happen that this element $g$ cannot be described 
in terms of the amalgam, as $G_{\pi_1}$ and $G_{\pi_2}$ might not be conjugate in the universal 
enveloping group of this amalgam, so that in this case it is very difficult to establish a nice correspondence between 
amalgams and coverings of geometries, as done in geometric covering theory. 
 If, however, $G_{\pi_2} = gG_{\pi_1}g^{-1}$ for some $g \in G_l$, then such a correspondence exists. In this case automatically $G_{l,\pi_2} = G_l \cap 
G_{\pi_2} = G_l \cap gG_{\pi_1}g^{-1} = gG_{l,\pi_1}g^{-1}$. Furthermore, $g \in G_l$ is an element of the 
amalgam, and we can {\em fuse} $G_{\pi_1}$ and $G_{\pi_2}$ via conjugation with $g$. We call such an amalgam 
$\mathcal{A}$ a {\em fused amalgam of parabolics} and depict it by 
$$ 
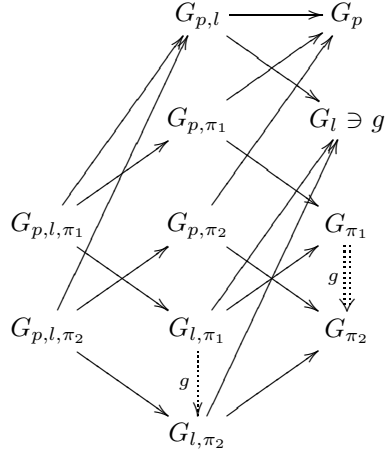
\xymatrix{ 
& G_{p,l} \ar[r] \ar[dr] & G_{p} \\ 
& G_{p,\pi_1} \ar[ur] \ar[dr]  & G_{l} \ni g \\ 
G_{p,l,\pi_1} \ar[uur] \ar[ur] \ar[dr] & G_{p,\pi_2} \ar[uur] \ar[dr] & G_{\pi_1} \ar@{:>}[d]_{g} \\ 
G_{p,l,\pi_2} \ar[uuur] \ar[ur] \ar[dr] & G_{l,\pi_1} \ar@{.>}[d]_{g} \ar[uur] \ar[ur] & G_{\pi_2} \\ 
& G_{l,\pi_2} \ar[ur] \ar[uuur] & } 
$$ 
 
The next definition formalises the concept of a fused amalgam. In this paper we only define fused amalgams sufficiently general for our purposes, although a number of possible generalisations come to mind immediately. 

\begin{definition}[Fused Amalgam] \label{fused}
Let ${\cal A} = ((G_j)_{j\in J}, (\phi_{ij})_{i<j})$ be an amalgam, with underlying graded poset ${\cal J} = (J, \leq)$ with grading function $\tau:J\rightarrow I = \{1,2,...,n\}$.
A {\em fusion} of $\cal A$, turning $\cal A$ into a {\em fused amalgam}, consists of three indices $j_0, j_1, j_2\in J_n = \tau^{-1}(n)$, an element $g \in G_{j_0}$, a lower neighbor $i_1$ of $j_0$ and $j_1$, a lower neighbor $i_2$ of $j_0$ and $j_2$, and an isomorphism $\gamma : G_{j_1} \to G_{j_2}$ such
that the following properties hold:
\begin{enumerate}
\item $\phi_{i_2,j_2}(gxg^{-1}) = (\gamma \circ \phi_{i_1,j_1})(x)$;
\item $\phi_{i,j_2}(x) =  (\gamma \circ \phi_{i,j_1})(x)$ for each $i < j_1, j_2$.
\end{enumerate}
\end{definition}

\begin{definition}[Enveloping Group] \label{completion}
Let $\A$ be a fused amalgam. A pair $(G,\pi)$ consisting of a group $G$ and a {map} $\pi: \sqcup\A\rightarrow 
G$ is called an {\em enveloping group} of $\A$, if 
\begin{enumerate} 
\item for all $j \in J$ the restriction of $\pi$ to $G_{j}$ is a homomorphism of $G_{j}$ to 
$G$; 
\item $\pi_{|G_{j}}\circ \phi_{i,j} = \pi_{|G_{i}}$ for all $i < j$; 
\item $\pi$ preserves fusion, i.e., $\pi(\gamma(x)) = \pi(g) \pi(x) \pi(g)^{-1}$ for every $x \in G_{j_1}$ and $g$, $\gamma$, $j_1$ as in Definition \ref{fused}; and 
\item $\pi(\sqcup \A)$ generates $G$. 
\end{enumerate} 
\end{definition}

\begin{proposition} \label{universal} 
Let $\A$ as above be a fused amalgam of groups, let $F(\A) = \gen{(u_g)_{g\in\A}}$ be the free group on the 
elements of $\A$ and let $$S_1 = \lb u_xu_y=u_z,\mbox{ whenever $xy=z$ in some $G_{j}$} \rb$$ and $$S_2 = \lb 
u_x = u_y,\mbox{ whenever $\phi(x)=y$ for some identification $\phi$} \rb $$ and $$S_3 = \lb u_x = 
gu_yg^{-1},\mbox{ whenever $x \in G_{j}$ and $y \in G_{j'}$ are fused by $g$} \rb $$ be relations for $F$. 
Then for each enveloping group $(G,\pi)$ of $\A$ there exists a unique group epimorphism $$\hat \pi : \mc{U}(\mc{A}) 
\rightarrow G$$ with  $\pi = \hat \pi \circ \psi$ where $$\mc{U}(\A) = \gen{(u_g)_{g\in\A} \mid S_1,S_2,S_3} 
\mbox{ and } \psi : \sqcup\mc{A} \rightarrow \mc{U}(\A) : g \mapsto u_g.$$ 
$$\xymatrix{ 
\sqcup\A \ar[r]^\psi \ar[dr]_\pi & \mc{U}(\mc{A}) \ar[d]^{\hat\pi} \\ 
&  G }$$ 
\end{proposition} 
 
\begin{proof} 
As in \cite{Gramlich/Maldeghem}. 
\end{proof} 
 
\begin{definition}[Universal Enveloping Group] 
Let $\A$ be a fused amalgam of groups. Then $$\psi : \sqcup\mc{A} \rightarrow \mc{U}(\A) : g \mapsto u_g$$ for 
$\mc{U}(\A)$ as in Proposition \ref{universal} is called the {\em universal enveloping group of $\A$}.
\end{definition}

\subsection{Some additional theory of intransitive geometries} 
 
Note that the covering theory from \cite{Gramlich/Maldeghem} does not apply to the geometry $\Gamma_2$ from 
Subsection \ref{gamma} by Proposition \ref{properties}. In this section we present a covering theorem making use of fused amalgams in order to 
tackle that geometry $\Gamma_2$.  For simplicity, we will state 
the theorem in such a way that it exactly fits the properties of $\Gamma_2$. Generalisations are of course 
possible. 
 
\begin{theorem} \label{mainthmcover} 
Let $\Gamma = (X,*,\typ)$ be a connected geometry over $I = \{ 1, 2, 3 \}$ having the direct sum property whose 
basic diagram $\Delta$ is $\node_1\arc\node_2\arc\node_3$. Let $G$ be a vertex-transitive group of 
automorphisms of $\G$ that acts transitively on the flags of type $\{ i, j \}$ for all $i, j \in I$ 
corresponding to adjacent nodes of the diagram $\Delta$. Furthermore, let $F = \{ w^1, w^2, w^3 \}$ be a flag, 
let $w^3$ and $gw^3$, $g \in G_{w^2}$, be orbit representatives of the action of $G_{w^1,w^2}$ on the elements 
of type $3$.
Finally, let $\mc{A} = \mc{A}(\Gamma,G,F)$ be the fused amalgam 
of parabolics. Then the diagram coset pregeometry $$\hat\Gamma = 
((\mc{U}(\mc{A})/G_{w^i}\times\{w^i\})_{i \in I},*)$$ is a simply connected geometry that admits a universal 
covering $\pi : \hat\Gamma \rightarrow \Gamma$ induced by the natural epimorphism $\mc{U}(\mc{A}) \rightarrow G$. 
Moreover, $\mc{U}(\mc{A})$ is of the form $\pi_1(\Gamma).G$. 
\end{theorem} 
 
\begin{proof} 
First notice that, since $\Gamma$ is connected, $G$ is generated by all its parabolics (different from $G$) by 
Theorem~\ref{char conn}. As the embedding of $\mc{A}$ in $G$ preserves fusion, by Definition~\ref{completion}, the group $G$ is an enveloping group of $\mc{A}$ and Proposition~\ref{universal} shows that the natural morphism $\mc{U}(\mc{A}) \rightarrow G$ is 
surjective. 
 
The map $$\phi : \sqcup\mc{A} \rightarrow G$$ and, thus, the map $$\hat\phi : \sqcup\mc{A} 
\rightarrow \mc{U}(\mc{A})$$ is injective. Therefore the natural epimorphism 
$$\psi: \mc{U}(\mc{A}) \rightarrow 
G$$ induces an isomorphism between the amalgam $\hat\phi(\sqcup\mc{A})$ inside $\mc{U}(\mc{A})$ and the amalgam 
$\phi(\sqcup\mc{A})$ inside $G$. Hence the epimorphism $\psi : \mc{U}(\mc{A}) \rightarrow G$ induces a quotient 
map between pregeometries $$\pi : \hat\Gamma = ((\mc{U}(\mc{A})/G_{w^i}\times\{w^i\})_{i \in I},*) \rightarrow 
((G/G_{w^i}\times\{w^i\})_{i \in I},*).$$ The latter diagram coset pregeometry is isomorphic to $\Gamma$ by the 
Reconstruction Theorem \ref{isomorphism theorem2}. Notice that $\mc{U}(\mc{A})$ acts on $\Gamma \cong 
((G/G_{w^i}\times\{w^i\})_{i \in I},*)$ via 
$$(gG_{w^i},w^i) \mapsto (\psi(u) gG_{w^i},w^i) \quad \mbox{ for } \quad u \in \mc{U}(\mc{A}).$$ 
 
We want to prove that this quotient map actually is a covering map. The pregeometry $\hat\Gamma$ is connected by 
Theorem \ref{char conn}, because $\mc{U}(\mc{A})$ is generated by $\hat{\phi}(\sqcup\mc{A})$. 
Let us start with proving the isomorphism between the residues of elements of type $2$.
Since the diagram is a string, a coset $xG_{w^1}$ is incident with a coset $yY$ with $Y = G_{w^3}$ or $Y = G_{gw^3}$, if and only if there exists a coset $zG_{w^2}$ such that $xG_{w^1}\cap zG_{w^2} \neq \emptyset \neq zG_{w^2}\cap yY$. So, all thats needs to be checked is that there is a bijection between the cosets of $G_{w^1}$ in ${\cal U}({\cal A})$ that meet $G_{w^2}$ and the corresponding cosets in $G$, and similarly for cosets of $Y$. However, the existence of such bijections is obvious. Indeed, the cosets of $G_{w^1}$ meeting $G_{w^2}$ are precisely those that can be written as $hG_{w^1}$ for $h\in G_{w^2}$, and similarly for cosets of $Y$.

Turning to the residue $w^1$, 
 let $hG_{w^3}$ represent a 3-element incident to $G_{w^1}$. Then there exists a coset $h_1G_{w^2}$ with $h_1\in G_{w^1}$ and $h_1G_{w^2}\cap hG_{w^3}\neq \emptyset$. As $h_1\in G_{w^1}$, we have
$G_{w^1}\cap h_1G_{w^2} = h_1G_{w^1,w^2}$. Turning to $hG_{w^3}$, the condition $h_1G_{w^2}\cap hG_{w^3}\neq\emptyset$ is equivalent to $G_{w^2}\cap h_1^{-1}hG_{w^3}\neq\emptyset$. This shows that we can choose $h$ such that $h_1^{-1}h = h_2\in G_{w^2}$, namely $h = h_1h_2\in G_{w^1}G_{w^2}$. Moreover, in view of the hypotheses we have assumed on $G_{w^1,w^2}$, there exists an element $f\in G_{w^1,w^2}$ such that either $fh_2G_{w^3} = G_{w^3}$ or $fh_2G_{w^3} = gG_{w^3}$. In the former case we can choose $h_2 = f^{-1}\in G_{w^1,w^2}$ and $h = h_1h_2\in G_{w^1}$. Thus, 
$G_{w^1}\cap hG_{w^3} = hG_{w^1,w^3}$. Also, $G_{w^1}\cap G_{w^2}\cap h_2G_{w^3}$ contains $h_2 f^{-1}\in G_{w^1,w^2}$. Hence $G_{w^1}\cap G_{w^2}\cap h_2G_{w^3} = f^{-1}G_{w^1,w^2,w^3}$. Accordingly, $G_{w^1}\cap h_1G_{w^2}\cap hG_{w^3} = h_1f^{-1}G_{w^1,w^2,w^3}$. So, these particular $\{2,3\}$-flags of $res(w^1)$ correspond to cosets of $G_{w^1,w^2,w^3}$ in $G_{w^1}$. However, the above holds in $\widehat{\cal G}$ as well as in $\cal G$. Consequently, that the part of $\widehat{\Gamma}_{w^1}$ formed by 2-elements and 3-elements in the same orbit of $G_{w^3}$ is isomorphic to the analogous part of ${\Gamma}_{w^1}$. 
Suppose now that the latter case occurs, namely $fh_2G_{w^3} = gG_{w^3}$. So, $hG_{w^3} = h_1h_2G_{w^3} = h_1f^{-1}gG_{w^3}$, whence $hG_{w^3}g^{-1} = h_1f^{-1}gG_{w^3}g^{-1}$. As $h_1G_{w^2} = h_1gG_{w^2}g^{-1}$ (because $gG_{w^2}g^{-1} = G_{w^2}$, since $g\in G_{w^2}$), we have $h_1G_{w^2}\cap hG_3\neq \emptyset$ if and only if $h_1G_2\cap h_1f^{-1}gG_{w^3}g^{-1}\neq\emptyset$. So, we can repeat the above argument with $hG_{w^3}$ replaced by $hG_{w^3}g^{-1} = h_1f^{-1}gG_{w^3}g^{-1}$, thus obtaining that the flag $\{h_1G_2, hG_3\}$ of $res(w^1)$ corresponds to a coset of $G_{w^1}\cap G_{w^2}\cap gG_{w^3}g^{-1}$ in $G_{w^1}$. In other words, the part of $res(w^1)$ formed by the 2-elements and the 3-elements of the orbit containing $gw^3$ is isomorphic to the geometry of cosets of $G_{w^1,w^2}$ and $G_{w^1,gw^3}$ inside $G_{w^1}$. Again, this is true in $\widehat{\Gamma}$ as well as in $\Gamma$. So, in either of these two geometries, that part of the residue of $w^1$ is canonically isomorphic to the same geometry of cosets inside $G_{w^1}$. So, that part of the residue of $w^1$ in $\widehat{\Gamma}$ is isomorphic to the corresponding part of the residue of $w^1$ in $\Gamma$.
So far, we have proved that each of the two parts of ${\widehat{\Gamma}}_{w^1}$ is isomorphic to the corresponding part in ${\Gamma}_{w^1}$. Moreover, it is clear from the above that the two `partial' isomorphisms constructed in this way from ${\widehat{\Gamma}}_{w^1}$ to ${\Gamma}_{w^1}$ agree on the set of 2-elements. Therefore they can be pasted together so that to construct an isomorphism from the whole of ${\widehat{\Gamma}}_{w^1}$ to the whole of ${\cal G}_{w^1}$. 

A similar argument applies to residues of $w^3$.  
 Hence $\pi : \hat\Gamma \rightarrow \Gamma$ induces isomorphisms between 
the residues of flags of rank one, so $\pi$ indeed is a covering of 
pregeometries. Since $\Gamma$ actually is a geometry the pregeometry $\hat\Gamma$ is also a geometry. 
 The universality of the covering $$\pi : \hat\Gamma \rightarrow \Gamma$$ induced by the canonical map $\mc{U}(\mc{A}) 
\rightarrow G$ is proved as in \cite[Theorem 3.1]{Gramlich/Maldeghem}. The structure of $\hat G \cong \mc{U}(\mc{A})$ is 
evident by combinatorial topology, cf.\ Chapter 8 of \cite{Seifert/Threlfall:1934}, restated in \cite[Section 2.2]{Gramlich/Maldeghem}. 
\end{proof} 
 
\begin{corollary}[Tits' lemma] \label{tits} \label{Tits2} \label{tits2} 
Let $\Gamma = (X,*,\typ)$ be a connected geometry over $I = \{ 1, 2, 3 \}$ having the direct sum property whose 
basic diagram $\Delta$ is $\node_1\arc\node_2\arc\node_3$. Let $G$ be a vertex-transitive group of 
automorphisms of $\Gamma$ that acts transitively on the flags of type $\{ i, j \}$ for all $i, j \in I$ 
corresponding to adjacent nodes of the diagram $\Delta$. Furthermore, let $F = \{ w^1, w^2, w^3 \}$ be a flag, 
let $w^3$ and $gw^3$, $g \in G_{w^2}$, be orbit representatives of the action of $G_{w^1,w^2}$ on the elements 
of type $3$.
Finally, let $\mc{A} = \mc{A}(\Gamma,G,F)$ be the fused amalgam of parabolics. 
 The geometry $\Gamma$ is 
simply connected if and only if the canonical epimorphism 
$$\mc{U}(\mc{A}(\Gamma,G,F)) \rightarrow G$$ is an 
isomorphism. \pend 
\end{corollary} 
 
\subsection{Amalgamation} \label{4.4}

By Proposition \ref{properties} the $\sfG_2(\mathbb{K})$ action on $\Gamma_2$ satisfies the hypotheses of Theorem \ref{mainthmcover}, so that by Proposition \ref{gamma3} we can apply Corollary \ref{tits} in order to obtain the following amalgamation result.
  
\begin{theorem} \label{amalg}
Let $p$, $l$, $\pi_1$, $\pi_2 = g\pi_1$, $g \in G_l$ be a set of orbit representatives of the $\sfG_2(\mathbb{K})$ action on $\Gamma_2$ such that $p$, $l$, $\pi_1$ is a chamber of $\Gamma_2$.
Then $\sfG_2(\mathbb{K})$ is the universal enveloping group of the following fused amalgam of parabolics:
$$ 
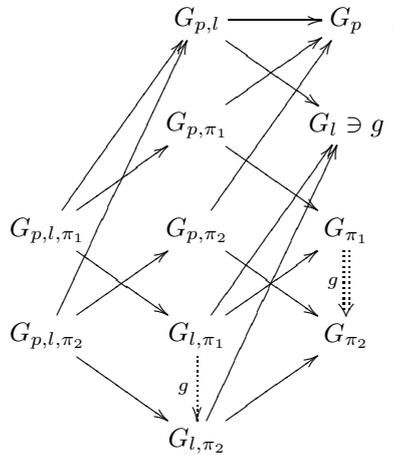
\xymatrix{ 
& G_{p,l} \ar[r] \ar[dr] & G_{p} \\ 
& G_{p,\pi_1} \ar[ur] \ar[dr]  & G_{l} \ni g \\ 
G_{p,l,\pi_1} \ar[uur] \ar[ur] \ar[dr] & G_{p,\pi_2} \ar[uur] \ar[dr] & G_{\pi_1} \ar@{:>}[d]_{g} \\ 
G_{p,l,\pi_2} \ar[uuur] \ar[ur] \ar[dr] & G_{l,\pi_1} \ar@{.>}[d]_{g} \ar[uur] \ar[ur] & G_{\pi_2} \\ 
& G_{l,\pi_2} \ar[ur] \ar[uuur] & }  .
$$ 
\end{theorem}

\section{Simple connectivity}\label{s5} 
 
In this section, we show the simply connectivity of most of the geometries $\Gamma_i$, $i=0,1,2,3$. To be precise, we prove that all these geometries are simply connected, except for $\Gamma_1$ with $|\K|=2$. 
 
Collinearity on $\sfQ(6,\K)$ and in $\mathsf{W}(5,\K)$ will be denoted by $\perp$. Also, in $\sfH(\K)$, there 
is a natural distance function on the set of points and lines, with values in $\{0,1,\ldots,6\}$ (this function 
is the graph theoretic distance in the incidence graph). Üairs of elements at distance $6$ will be called 
\emph{opposite}. In $\PG(5,\K)$, there are two types of singular planes: the ideal planes, and the hexagonal 
singular planes. If $\alpha$ is a hexagonal singular plane, then the ideal centres of all ideal lines in 
$\alpha$ coincide and will be called the \emph{hexagonal pole} of $\alpha$ (it is the intersection of all 
hexagon lines in $\alpha$). The pencil of hexagon lines will be denoted by $\mathcal{P}_\alpha$. Also, a 
nonsingular plane $\beta$ contains a pencil of symplectic lines; this pencil will be denoted by 
$\mathcal{P}_\beta$ and the intersection of the lines of $P_\beta$ shall be called the \emph{pole} of $\beta$. 
 
 
First we will start by proving simple connectivity of $\Gamma_0$, $\Gamma_1$, $\Gamma_2$ as defined in 
Subsection~\ref{gamma}. In what follows, we will use over and over the simple observation that, if two 
symplectic lines $L_1,L_2$ meet in a point $p$, then every line in the plane $\langle L_1,L_2\rangle$ through 
$p$ is symplectic. Moreover, if the plane $\langle L_1,L_2\rangle$ is nondegenerate, then at most one of these 
lines can be hexagonal, because a pair of intersecting hexagon lines spans a totally isotropic plane, see Section \ref{2.1}, item (i).
 In the following, a \emph{geometric triangle} is a triangle consisting of points and lines in a plane of the 
geometry. 

The next proposition has also been proved by Blok and Hoffman \cite{Blok/Hoffman}, but we provide a short proof for 
sake of completeness.

\begin{proposition} \label{gamma1}
The geometry $\Gamma_0$ is simply 
connected. 
\end{proposition} 
 
\begin{proof} 
Clearly the diameter at type 1 of the truncation of $\Gamma_0$ at 
the types 1 and 2 is two, so it suffices to show that we can 
subdivide any triangle, any quadrangle and any pentagon in geometric triangles. 
Every triangle is geometric, so for triangles there is nothing to prove. 
 
Let $a,b,c,d$ be a quadrangle of type 1 elements (points of $\PG(5,q)$). We may assume that $ac$ and $bd$ are 
symplectic lines, since otherwise the quadrangle automatically decomposes into triangles. Choose any point $e$ 
on $ab$. If $ce$ were symplectic, then so would $cb$. Hence $ce$, and similarly also $de$, are 
non-symplectic. We have subdivided our quadrangle in the triangles $a,d,e$ and $c,d,e$ and $b,c,e$. 
 
Now let $a,b,c,d,e$ be a pentagon of type 1 elements. Again we may 
assume that all of $ac$, $bd$, $ce$, $da$ and $eb$ are symplectic, since otherwise the pentagon decomposes automatically. 
Since $a,c,d$ is a triangle in $\PG(5,q)$, the corresponding 
symplectic hyperplanes $a^\perp$, $c^\perp$, $d^\perp$ do not meet in a $3$-space, whence their union cannot cover the whole space. Therefore there is 
a point $f$ with $cf$, $df$ and $af$ non-symplectic lines and we have 
subdivided our pentagon into the null-homotopic circuits 
$a,b,c,f$ and $c,d,f$ and $d,e,a,f$. 
\end{proof} 
 
\begin{proposition} \label{gamma2}
The geometry $\Gamma_1$ is simply connected, 
whenever $|\K|>2$. 
\end{proposition} 
 
\begin{proof} 
Since the point-line truncations of $\Gamma_0$ and $\Gamma_1$ coincide, by the proof of 
Proposition~\ref{gamma1}, it suffices to show that every triangle is null-homotopic. 
 
Let $a,b,c$ be a triangle, and suppose it is not geometric. Hence $a,b,c$ are three pairwise opposite points in 
$\sfH(\K)$ and the plane $\< a,b,c\>$ contains a (unique) hexagon line $l$. We may assume that $a$ is not 
incident with $l$. Then there exists a hexagon line $m$ through $a$ not concurrent with $l$ and not concurrent 
with a hexagon line that is concurrent with $l$. It follows that the $3$-space $\Xi:=\<l,m\>$ is nondegenerate 
and contains a regulus, consisting of hexagon lines, of a ruled nondegenerate quadric $Q$, cf.\ Section \ref{2.1}, item (v). Denote the 
unique line of $Q$ in $\<a,b,c\>$ different from $l$ by $l'$ (and note that $l'$ is incident with $a$ and meets $l$), and for 
each point $z\in l'$, denote by $l_z$ the unique hexagon line on $Q$. To prove the claim, it suffices to find a 
point $x$ such that the planes $\<a,b,x\>,\<a,c,x\>$ and $\<b,c,x\>$ do not contain any hexagon line and such 
that the lines $ax,bx,cx$ are not symplectic. The latter is satisfied whenever $x \in \Xi$ is not contained in the 
union $\mathcal{U}_1$ of the planes $\pi_a,\pi_b,\pi_c$, where $\pi_a$ is generated by the symplectic lines 
through $a$ inside $\Xi$, and likewise for $\pi_b$ and $\pi_c$. Note that $\gen{a,b,c} = \gen{l,l'} \subseteq \Xi$. The former is satisfied whenever $x$ does 
not lie in the union $\mathcal{U}_2$ of the planes $\<a,b,l_{ab\cap l'}\>,\<a,c,l_{ac\cap l'}\>,\<b,c,l_{bc\cap 
l'}\>$ and $\<a,b,c\>$. Since $l'$ contains $a$, we see that $l_{ab\cap l'}=l_{ac\cap l'}=m$. 
 Since $\Xi$ cannot be the union of seven planes if $|\K| \geq 8$, we 
may suppose $|\K|=4$. In that case 
$\mathcal{U}_2$, which is the union of four planes no three of which meet in a 
line, i.e., a tetrahedron, covers exactly 58 points, 
leaving a set $S$ of $85-58 = 27$ possibilities for $x$. If a plane contained in 
$\mathcal{U}_1$ does not contain an 
intersection line of two planes in $\mathcal{U}_2$, then it meets $S$ in $7$ 
or $6$ points. If, on the other 
hand, such a plane does contain such an intersection line, it contains $9$ 
points of $S$. Consequently we may 
assume that the planes in $\mathcal{U}_1$ partition $S$ and each plane 
contains some intersection line of 
planes in $\mathcal{U}_2$. For all three planes, 
this line must be the same, as otherwise the three intersections would have to be pairwise distinct. Hence in that case each of the three planes would have to contain two lines of the tetrahedron, which would imply that they actually belong to $\mathcal{U}_2$, a contradiction. But 
then the planes $\pi_a$, $\pi_b$, $\pi_c$ contain a common line, which implies 
that $a,b,c$ are collinear in 
$\PG(5,\K)$ (since $\Xi$ is nondegenerate), another contradiction. 
\end{proof} 
 
The above proof fails for $|\K|=2$. In fact, a computation using \textsf{GAP} reveals that $\Gamma_1$ admits in 
this case a $3$-fold universal cover. The source code of the program we used for verification of this fact can be found in 
\ifthenelse{\boolean{isArxivVersion}}{Appendix \ref{appendix:gap-code} and also on the website \cite{Horn}}{\cite{arxiv}}.
 
\begin{proposition} \label{gamma3}
The geometry $\Gamma_2$ is simply connected. 
\end{proposition} 
  
\begin{proof}
As in the case of $\Gamma_1$ we only need to prove that every triangle is null-homotopic. 
Suppose $|\K| > 2$ (hence $|\K| \geq 4$, as the characteristic of $\K$ is two). Given a triangle 
$p_0,p_1,p_2$ of the collinearity graph of $\Gamma_2$, not contained in one line of $\Gamma_2$, set $\pi := 
\langle p_0, p_1, p_2\rangle$. Clearly, the plane $\pi$ is nondegenerate. Let $p$ be its pole. 
 We may assume that $\pi$ is a plane of $\Gamma_1$, as otherwise our triangle is geometric. So all the 
symplectic lines of $\pi$ are ideal. Given a line $l\in{\cal P}_\pi$, let $\pi_l$ be the hexagonal 
plane on $l$ and let $p_l = p_{\pi_l}$ be the hexagonal pole of $\pi_l$. 
 Suppose, by way of contradiction, that $\pi\subset p_l^\perp$. So, the singular planes on $pp_l$ are precisely 
those spanned by $p_l$ and a line $m\in {\cal P}_\pi$. Since $pp_l$ is a hexagon line, all of these planes are 
hexagonal and the map sending a line $m\in{\cal P}_\pi$ to the hexagonal pole $p_m$ of $\langle m,p_l\rangle$ 
is a bijection from ${\cal P}_\pi$ to the set of points of $pp_l$. Therefore, $p = p_m$ for some line 
$m\in{\cal P}_\pi$. Hence $m$ is a hexagon line, contrary to our assumptions. 
 As a consequence, $p_l^\perp\cap \pi = l$ for every line $l\in{\cal P}_\pi$. We can now choose the line $l$ 
such that $l\in{\cal P}_\pi\setminus\{pp_0, pp_1, pp_2\}$ (noting that we have assumed $|\K| > 2$). For $1\leq 
i<j\leq 3$, none of the planes $\pi_{ij} = \langle p_i, p_j, p_L\rangle$ is singular (since they contain the 
non-symplectic line $p_ip_j$), but each of them contains a hexagon line, namely the line $a_{ij}p_l$, where 
$a_{ij} = l\cap p_ip_j$. 
 So, we have decomposed $p_0,p_1,p_2$ into triangles $p_i,p_j,a$, each of which is contained in a plane of 
$\Gamma_2$. 
 
For $|\K|=2$, a computer based argument proves the claim. Again, the source code for the \textsf{GAP} program we used can be found in 
\ifthenelse{\boolean{isArxivVersion}}{Appendix \ref{appendix:gap-code} and also on the website \cite{Horn}}{\cite{arxiv}}.
\end{proof}

\begin{proposition}\label{simpconn:Gamma3} \label{prop3}
The geometry $\Gamma_3$ is simply connected. 
\end{proposition} 
 
\begin{proof} We prove this in a series of lemmas. The strategy is to show that every cycle in the collinearity 
graph of $\Gamma_3$ is null homotopic. We begin by noting that the diameter of that graph is equal to two. 
Indeed, if two points $a,b$ are incident with the same hexagon line $l$, then we can choose a hexagonal plane 
$\pi$ through $l$ such that neither $a$ nor $b$ is incident with at least two hexagon lines of $\pi$. Consequently, for 
any point $c$ in $\pi$ not on $l$ the lines $ac$ and $bc$ are ideal. If two points $a,b$ are at distance two in 
the collinearity graph of $\sfQ(6,\K)$, then the points $c$ collinear to both $a$ and $b$ form a generalised 
quadrangle $\sfQ(4,\K)$; those for which either $ac$ or $bc$ are hexagon lines form two lines in $\sfQ(4,\K)$. 
Hence there are plenty of points $c$ in $\sfQ(4,\K)$ for which $ac$ and $bc$ are ideal. 
 As a consequence, we only have to show that triangles, quadrangles and pentagons are null homotopic. This will be done in lemmas \ref{triangles}, \ref{quadranglesBIS}, and \ref{pentagons}. 
\end{proof}

\begin{lemma} Every triangle is null homotopic.\label{triangles}
\end{lemma} 
 
\begin{proof} 
If the triangle is geometric, i.e., contained in an ideal plane, then this 
is trivial. So we may assume that 
we have a triangle $a,b,c$ in a hexagonal plane.
Let $\pi'$ be an ideal plane on some ideal line of $\gen{a,b,c}$ not containing $a$, $b$ or $c$. Then, by Section \ref{2.1}, item (iv), the ideal centres of the ideal lines of $\pi'$ form an ideal plane $\pi$. Then span $H := \gen{\pi,\pi'}$ is a hyperplane of $\PG(6,\K)$ meeting $Q$ in a nondegenerate 
hyperbolic quadric $Q^+$. Moreover, none of the points $a$, $b$, $c$ is incident with $\pi \cup \pi'$, and $\pi$ contains the pole of the plane $\gen{a,b,c}$. By Section \ref{2.1}, item (iv) the hexagonal planes in $H$ are those that share a point 
of $\pi$ or $\pi'$ and a line of 
$\pi'$ or $\pi$, respectively.

Now we apply the Klein correspondence to view the situation in a $3$-space 
$\PG(3,\K)$. To be precise, we map $\pi$ to a point $x$ and hence $\pi'$ to a plane $\alpha$ off that point. Translated to the space $\PG(3,K)$, the hexagonal planes in $H$ correspond with the points of $\alpha$ and 
with the planes through $x$. 
The plane $\gen{a,b,c}$ corresponds to a point $z$ in $\alpha$, and 
the points $a,b,c$ correspond to 
lines $A,B,C$, respectively, through $z$, but not contained in $\alpha$ and 
not incident with $x$. Also, none 
of the planes $\gen{A,B}$, $\gen{A,C}$, $\gen{B,C}$ contain $x$. Let $l$ be 
a line in $\alpha$ not through $z$ 
and let $\beta$ be a plane through $l$, different from $\alpha$, and not 
incident with $x$. Then the 
intersections $L_C:=\beta\cap \gen{A,B}$, $L_B:=\beta\cap\gen{A,C}$ and 
$L_A:=\beta\cap\gen{B,C}$ form a 
triangle which corresponds with a null-homotopic triangle in $\Gamma$. By 
construction also the triangles 
$A,B,L_C$ and $A,C,L_B$ and $B,C,L_A$ correspond with null-homotopic triangles 
in $\Gamma$, and likewise so do the 
triangles $L_A,L_B,C$ and $L_A,L_C,B$ and $L_B,L_C,A$. Hence also the triangle 
$a,b,c$ is null homotopic as it it 
is the sum of seven geometric triangles. 
\end{proof} 
 
\begin{lemma} Every quadrangle $a,b,c,d$ with $a\not\perp 
c$ and $b\not\perp d$, is null homotopic.\label{quadrangles} 
\end{lemma} 
 
\begin{proof} 
Let $x$ be the ideal centre of $ab$ and let $y$ be the ideal centre of $cd$. It is easy to see that our 
assumptions imply that $x$ and $y$ are not collinear in $\sfH(\K)$. Suppose now first that the planes 
$\gen{x,a,b}$ and $\gen{y,c,d}$ are disjoint. This, as $x$, $y$ are noncollinear in $\sfH(\K)$ is equivalent to $x$ being opposite $y$ in $\sfH(\K)$. Let $X_1$ be the set of points of $\sfH(\K)$ collinear with $x$ 
and not opposite $y$, and let likewise $Y_1$ be the set of points of $\sfH(\K)$ not opposite $x$ and collinear 
with $y$. Note that $a\in X_1$ if and only if $d\notin Y_1$ (and similarly, $b\in X_1$ if and only if $c\notin 
Y_1$). For assume $a\in X_1$, then clearly, the only point of $Y_1$ collinear in $\sfQ(6,\K)$ with $a$ is also 
collinear with $a$ in $\sfH(\K)$, thus not collinear with $a$ in $\Gamma$, and hence $d\notin Y_1$. The 
remaining implications follow identically. 
In view of the above, either $b \not\in X_1$ or $c \not\in Y_1$. We may assume that $b \not\in X_1$. So $c \in Y_1$ and we are left with the following cases.

\begin{enumerate} 
\item In case $a \not\in X_1$, the intersection $a^\perp\cap b^\perp\cap\gen{y,c,d}$ is a singleton 
$\{u\}$. Our assumptions imply that $u\notin\{c,d,y\}$. If $u$ belonged to the (hexagon) line $yc$, then, 
since $c\in b^\perp$, also $y\in b^\perp$, a contradiction. So $u\notin yc$ and likewise $u\notin yd$. Moreover, neither $au$ nor $bu$ are hexagon lines, as neither $a$ nor $b$ belong to $X_1$. Hence 
$u$ is collinear in $\Gamma$ with all of $a,b,c,d$ and we can subdivide the quadrangle $a,b,c,d$ in the triangles 
$a,b,u$ and $b,c,u$ and $c,d,u$ and $d,a,u$. 
 \item In case $a\in X_1$, we have $d\notin Y_1$. Pick any point $v$ 
on the ideal line $ab$, different from $a$ and $b$. Our assumptions imply easily that 
there is a unique point $w$ on $cd$ collinear in $\Gamma$ with $v$. Since $v\notin X_1$, we know that $vw$ 
is an ideal line. Also, since $c$ is collinear with $b$ in $\Gamma$, and $w$ cannot be, we deduce that $w\neq c$. 
Similarly $w\neq d$. By the previous arguments the quadrangles $a,v,w,d$ and $v,b,c,w$ are null homotopic, and they 
subdivide $a,b,c,d$. Hence also in this case $a,b,c,d$ is null homotopic. Note that the situation considered here does not occur when $|\K|=2$. 
 \end{enumerate} 
 
Suppose now secondly that the planes $\gen{a,b,x}$ and $\gen{c,d,y}$ meet in a point $z$. Our assumptions imply that 
$z$ is distinct from the intersection point $x_1$ of $ab$ with $xz$, and also from the intersection point $y_1$ 
of $cd$ with $yz$. Now clearly the ideal centres of $ad$ and $bc$ are opposite $z$, which is the ideal centre of $x_1y_1$. 
Hence, by the previous arguments (with $ab$ replaced by $wd$ and $vb$, respectively), the quadrangles $a,x_1,y_1,d$ and $b,x_1,y_1,d$ are null homotopic. Since they subdivide $a,b,c,d$, the lemma follows. 
 \end{proof} 
 
\begin{lemma} Every quadrangle $a,b,c,d$ is null homotopic.\label{quadranglesBIS} 
\end{lemma} 
 
\begin{proof} 
By Lemma \ref{quadrangles}, we may assume that $a$ and $c$ are 
collinear on $\sfQ(6,\K)$. If $ac$ is ideal, then we are done by the 
fact that all triangles are null homotopic, see Lemma 
\ref{triangles}. Hence we may assume that $ac$ is a hexagon line. 
Clearly, we may assume that $b$ and $d$ are not collinear on the 
quadric as otherwise $a,b,c,d$ lie in a plane of the quadric and 
then $ad$ meets $bc$ in some point $e$. The triangles $a,b,e$ and 
$c,d,e$ are null homotopic by Lemma \ref{triangles}, hence the 
result. 
 
The plane $\gen{a,b,c}$ is degenerate since it contains a hexagon 
line. Considering any other plane $\pi$ through $bc$, it follows 
that $\pi$ is an ideal plane. Since $d\not\perp b$, the line 
$l=\pi\cap d^\perp$ does not coincide with $bc$ and hence contains 
at least two points $u,v$ off $bc$. If both $du$ and $dv$ were 
hexagon lines, then also $dc$ must be a hexagon line, a 
contradiction. So we may assume that $du$ is ideal. But now we 
have subdivided the quadrangle $a,b,c,d$ into the circuits 
$c,d,u$ and $c,u,b$ and $u,b,a,d$. The two former are null homotopic 
by Lemma~\ref{triangles}, while the latter is null homotopic by 
Lemma~\ref{quadrangles}, noting that $u$ is not collinear with $a$ 
on the quadric. 
\end{proof}

\begin{lemma} 
Every pentagon $a_0,a_1,a_2,a_3,a_4$ is null homotopic. \label{pentagons} 
\end{lemma} 
 
\begin{proof} 
For $i = 0, 1,..., 4$, let $l_i = a_ia_{i+1}$ be the ideal line through $a_i$ and $a_{i+1}$ (indices being computed modulo 5). Suppose first that $l_{i+2}\not\subset a_i^\perp$ for some $i$. For instance, $l_2\not\subset a_0^\perp$. As the line $l_2$ is ideal, all singular planes on $l_2$ but one are ideal. So, we can choose an ideal plane $\beta$ on $l_2$ such that the plane $\alpha = \langle a_0^\perp\cap\beta, a_0\rangle$ is different from the unique hexagonal plane containing all hexagonal lines through $a_0$. Consequently, at most one of the lines of $\alpha$ through $a_0$ is hexagonal. Given an ideal line $l$ of $\alpha$ through $a_0$, let $c = l\cap(\alpha\cap\beta)$. Then $c$ is collinear with both $a_2$ and $a_3$ in $\Gamma_3$, since $\beta$ is an ideal plane. Thus, we can split $a_0, a_1, a_2, a_3, a_4$ into $a_0, a_1, a_2, c$ and $c, a_2, a_3$ and $a_0, c, a_3, a_4$. 
 
Suppose now that $l_{i+2}\subset a_i^\perp$ for every $i = 0, 1,..., 4$. Then $a_0, a_1, a_2, a_3, a_4$ is contained in a singular plane. Put $b := l_0\cap l_2$ and $c := l_4\cap l_2$. Then $a_0, a_1, a_2, a_3, a_4$ splits into the sum of $a_0, b, c$ and $b,a_1,a_2,b$ and $c,a_3,a_4$. 
 \end{proof} 
 
\ifthenelse{\boolean{isArxivVersion}}{

\appendix

\section{GAP code} \label{appendix:gap-code}

\begin{lstlisting}[frame=single,caption=gamma1.gap]
# Setup
Read("util.gap");
F := GF(q);
z := Z(q);
G2q := G2(q);
v_orb := [0,0,1,0,1,0] * z^0;
Display("Verifying group size...");
Assert(0, Size(G2q) = q^6*(q^6-1)*(q^2-1));

# The second (bigger) geometry, with the second plane orbit,
# is determined by the following elements:
u_space := [ 
		[[1,0,0,0,0,0]],
		[[1,0,0,0,0,0],                [0,0,0,1,0,0]],
		[[1,0,0,0,0,0], [0,0,1,0,0,1], [0,0,0,1,0,0]]
	] *z^0;
u_stab := [];
u_index := [
	(q^6-1)/(q-1),
	q^4*(1+q^2+q^4),
	(1+q)*(1+q^2+q^4)*(q^4-q^3)
	];
u_size := List(u_index, x -> Size(G2q) / x);

Read("build_stabs.gap");
Read("build_amalgam2.gap");

Display("Coset enumeration...");
A := FG / rels;
U := Subgroup(A, GeneratorsOfGroup(A){stab_gens[1]});
idx := Index(A, U);

Print("Expected index = ", Index(G2q, U1), "\n");
Print("Actual index = ", idx, "\n");
Print("Resulting covering factor = ", idx/Index(G2q, U1), "\n");
\end{lstlisting}
\begin{lstlisting}[frame=single,caption=gamma2.gap]
# Setup
Read("util.gap");
F := GF(q);
z := Z(q);
G2q := G2(q);
v_orb := [0,0,1,0,1,1] * z^0;
Display("Verifying group size...");
Assert(0, Size(G2q) = q^6*(q^6-1)*(q^2-1));

# The first (smaller) geometry, with the first plane orbit, 
# is determined by the following elements:
u_space := [ 
		[[1,0,0,0,0,0]],
		[[1,0,0,0,0,0],                [0,0,0,1,0,0]],
		[[1,0,0,0,0,0], [0,1,0,0,0,0], [0,0,0,1,0,0]],
		[[1,0,0,0,0,0],                [0,0,0,1,0,0], [0,0,0,0,0,1]]
	] * z^0;
u_stab := [];
u_index := [
	(q^6-1)/(q-1),
	q^4*(1+q^2+q^4),
	(1+q)*(1+q^2+q^4)*(q^3+q^2),
	(1+q)*(1+q^2+q^4)*(q^3+q^2)
	];
u_size := List(u_index, x -> Size(G2q) / x);

Read("build_stabs.gap");
Read("build_amalgam.gap");

Display("Coset enumeration...");
A := FG / rels;
U := Subgroup(A, GeneratorsOfGroup(A){stab_gens[1]});
idx := Index(A, U);

Print("Expected index = ", Index(G2q, U1), "\n");
Print("Actual index = ", idx, "\n");
Print("Resulting covering factor = ", idx/Index(G2q, U1), "\n");
\end{lstlisting}
\begin{lstlisting}[frame=single,caption=build\_stabs.gap]
for i in [1..3] do
	Print("Computing stabilizer number ", i, "\n");
	if q=2 then
		u_stab[i] := Stabilizer(G2q, u_space[i], OnSubspacesByCanonicalBasis);
	else
		# Computing the full stabilizer of u_space[i] directly is very expensive.
		# So instead, we use a trick: We compute the stabilizer of u_space[i]
		# inside two small subgroups, and then take the group generated
		# by these two stabilizers. Empirically, this gives us the desired
		# full stabilizer. We verify that we obtained the correct subgroup
		# by checking the size of the resulting group.
		if q = 4 then
			tmp1 := Group(G2q.1, G2q.2);
		else
			tmp1 := Group(G2q.1, G2q.3);
		fi;
		tmp2 := Group(G2q.2, G2q.3);
		tmp1 := Stabilizer(tmp1, u_space[i], OnSubspacesByCanonicalBasis);
		tmp2 := Stabilizer(tmp2, u_space[i], OnSubspacesByCanonicalBasis);
	
		u_stab[i] := Group(Union(List([tmp1, tmp2], GeneratorsOfGroup)));
		sx := SizeViaOrbit(u_stab[i], v_orb);
		Assert(0, sx = u_size[i]);
		SetSize(u_stab[i], u_size[i]);
	fi;
od;
U1 := u_stab[1];
U2 := u_stab[2];
U3 := u_stab[3];

if Size(u_size) > 3 then
	Display("Computing the second plane stabilizer");
	# We know that U3 and U4 are conjugated by the following element
	# in the line stabilizer:
	g := [
		[0,0,0,1,0,0],
		[0,0,0,0,0,1],
		[0,0,0,0,1,0],
		[1,0,0,0,0,0],
		[0,0,1,0,0,0],
		[0,1,0,0,0,0]]*z^0;

	U4 := U3^g;
fi;
\end{lstlisting}
\begin{lstlisting}[frame=single,caption=build\_amalgam.gap]
# U3 and U4 are conjugated by this element in U2 (which
# exchanges the two plane classes, and preserves lines).
g := [
	[0,0,0,1,0,0],
	[0,0,0,0,0,1],
	[0,0,0,0,1,0],
	[1,0,0,0,0,0],
	[0,0,1,0,0,0],
	[0,1,0,0,0,0]]*z^0;

#
# Compute the stabilizer intersections
#
Display("Computing the double stabilizers...");
U12 := Stabilizer(U2, u_space[1], OnSubspacesByCanonicalBasis);

U13 := Stabilizer(U3, u_space[1], OnSubspacesByCanonicalBasis);
U23 := Stabilizer(U3, u_space[2], OnSubspacesByCanonicalBasis);
U123 := Stabilizer(U23, u_space[1], OnSubspacesByCanonicalBasis);

U14 := Stabilizer(U4, u_space[1], OnSubspacesByCanonicalBasis);
U24 := Stabilizer(U4, u_space[2], OnSubspacesByCanonicalBasis);
U124 := Stabilizer(U24, u_space[1], OnSubspacesByCanonicalBasis);


Display("Finding generators suitable for amalgamation...");

repeat
	x := PseudoRandom(U12);
until ClosureGroup(U123,x)=U12;

repeat
	y := PseudoRandom(U23);
until ClosureGroup(U123,y)=U23;

repeat
	a := PseudoRandom(U13);
until ClosureGroup(U123,a)=U13 and ClosureGroup(U23,a)=U3;

repeat
	c := PseudoRandom(U14);
until ClosureGroup(U124,c)=U14;

repeat
	h := PseudoRandom(U124)^(g^-1);
until U124 = ClosureGroup(U123^g,h^g);


# Verify "nice" generator properties -- namely, that the intersections of
# the groups are generated by the intersections of their generating sets.
Display("Verifying...");

Assert(0, U12 = ClosureGroup(U123,x));
Assert(0, U23 = ClosureGroup(U123,y));

Assert(0, U12 = ClosureGroup(U124,x));
Assert(0, U24 = ClosureGroup(U124,y^g));

Assert(0, Size(U1) = SizeViaOrbit(ClosureGroup(U123,[x,a]),v_orb));
Assert(0, U2 = ClosureGroup(U123,[x,y,y^g]));
Assert(0, U3 = ClosureGroup(U123,[y,a]));
Assert(0, U4 = ClosureGroup(U124,[y^g,a^g]));


Display("Building the amalgam...");

g123 := ShallowCopy(SmallGeneratingSet(U123));
g123_names := List([1..Length(g123)], n -> Concatenation("f", String(n)));
g124_names := List(g123_names, x -> Concatenation(x, "G"));

gens := Concatenation(
			g123,
			List(g123, e->e^g),
			[h,h^g,x,y,y^g,a,a^g,c,c^g,g]);
gen_names := Concatenation(
			g123_names,
			g124_names,
			["h", "hG", "x", "y", "yG", "a", "aG", "c", "cG", "g"]);

Assert(0, Size(gens) = Size(gen_names));

u_gens := [];
u_gens[1] := Concatenation(g123_names, g124_names, ["hG","x","a", "c"]);
u_gens[2] := Concatenation(g123_names, g124_names, ["h","hG","x","y","yG","g"]);
u_gens[3] := Concatenation(g123_names,             ["h","y","a", "cG"]);
u_gens[4] := Concatenation(            g124_names, ["hG","yG","aG", "c"]);

# We provided lists of generators identified with elements in G2. Now
# we compute all relators between these.
Read("build_rels.gap");

# Identify generators with their conjugates (h/hG, y/yG, a/aG, c/cG), by adding
# suitable relators. For this it is sufficient to add the relations between
# a / a^g, and c / c^g, since the others are contained in U2 together with g,
# so the corresponding relators are already present.

fg_gens := GeneratorsOfGroup(FG);
g_in_FG := fg_gens[Position(gens, g)];
for e in ["a", "c"] do
	tmp1 := fg_gens[Position(gen_names, e)];
	tmp2 := fg_gens[Position(gen_names, Concatenation(e, "G"))];
	Add(rels, tmp2^-1 * tmp1^g_in_FG );
od;
\end{lstlisting}
\begin{lstlisting}[frame=single,caption=build\_amalgam2.gap]
#
# Compute the stabilizer intersections
#
Display("Computing the double stabilizers...");
U12 := Stabilizer(U2, u_space[1], OnSubspacesByCanonicalBasis);
U13 := Stabilizer(U3, u_space[1], OnSubspacesByCanonicalBasis);
U23 := Stabilizer(U3, u_space[2], OnSubspacesByCanonicalBasis);
U123 := Stabilizer(U23, u_space[1], OnSubspacesByCanonicalBasis);


Display("Finding generators suitable for amalgamation...");

# Find x in U12 such that U12 = <U123, x>
repeat
	x := PseudoRandom(U12);
until ClosureGroup(U123,x)=U12;

# Find y in U23 such that U23 = <U123, y>
repeat
	y := PseudoRandom(U23);
until ClosureGroup(U123,y)=U23;

# Find r in U13 such that U1 = <U12, r> and U3 = <U23, r>
repeat
	r := PseudoRandom(U13);
until ClosureGroup(U12,r)=U1 and ClosureGroup(U23,r)=U3;


# Verify "nice" generator properties -- namely, that the intersections of
# the groups are generated by the intersections of their generating sets.
Display("Verifying...");
Assert(0, U12 = ClosureGroup(U123,x));
Assert(0, U23 = ClosureGroup(U123,y));
Assert(0, U13 = ClosureGroup(U123,r));

Assert(0, U1 = ClosureGroup(U123,[x,r]));
Assert(0, U2 = ClosureGroup(U123,[x,y]));
Assert(0, U3 = ClosureGroup(U123,[y,r]));

Assert(0, Size(U1) = SizeViaOrbit(U1, v_orb));
Assert(0, Size(U2) = SizeViaOrbit(U2, v_orb));
Assert(0, Size(U3) = SizeViaOrbit(U3, v_orb));


Display("Building the amalgam...");

g123 := ShallowCopy(SmallGeneratingSet(U123));
g123_names := List([1..Length(g123)], n -> Concatenation("f", String(n)));

gens := Concatenation(g123, [x,y,r]);
gen_names := Concatenation(g123_names, ["x", "y", "r"]);

Assert(0, Size(gens) = Size(gen_names));

u_gens := [];
u_gens[1] := Concatenation(g123_names, ["x","r"]);
u_gens[2] := Concatenation(g123_names, ["x","y"]);
u_gens[3] := Concatenation(g123_names, ["y","r"]);

# We provided lists of generators identified with elements in G2. Now
# we compute all relators between these.
Read("build_rels.gap");
\end{lstlisting}
\begin{lstlisting}[frame=single,caption=build\_rels.gap]
Print("Constructing the amalgam...\n");

SetInfoLevel(InfoFpGroup, 3);

#FG := FreeGroup(Length(gens));
Assert(0, Length(gens) = Length(gen_names));
FG := FreeGroup(gen_names);

# Select the generators for each stabilizer
stab_gens := List(u_gens, gs -> List(gs, g -> Position(gen_names, g)));


u_rels := [];
rels := [];
stabs := List( stab_gens, t->Group( gens{t} ) );

# We now determine finite presentations for each stabilizer.
# To do this, we first find a permutation group isomorphic
# to that stabilizer, then use IsomorphismFpGroupByGenerators
# to get the finite presentation.
#
# Finally, we form the union of the relators of the stabilizers.
#
for i in [1..Length(stab_gens)] do
	# Verify and set the size of the group.
	# We assume here that v_orb has been set to a vector with short orbit!
	Assert(0, SizeViaOrbit(stabs[i], v_orb) = u_size[i]);
	SetSize(stabs[i], u_size[i]);
	#
	Print("Generating relators for ", stab_gens[i], "\n");
	#
	# Determine a nice permutation presentation
	# We assume here that v_orb has been set to a vector with short orbit!
	phi := IsomorphismPermGroupViaOrbit(stabs[i], v_orb);
	H := Image(phi);
	tmp_FG_gens := GeneratorsOfGroup(FG){stab_gens[i]};
	u_rels[i] := GroupRelatorsViaOrbit(H, GeneratorsOfGroup(H), tmp_FG_gens);
	# Verify all relations actually hold in G
	Assert(0, ForAll( u_rels[i], r->MappedWord( r, tmp_FG_gens,
						gens{stab_gens[i]} ) = One(stabs[i]) ));
	#
	Append(rels, u_rels[i]);
od;

rels := Set(rels);

Assert(0, ForAll( rels, r->MappedWord( r, GeneratorsOfGroup(FG), gens ) = gens[1]^0 ));
\end{lstlisting}
\begin{lstlisting}[frame=single,caption=util.gap]
#
# We can efficiently compute a (lower bound on) the size of a matrix group G
# by finding a vector v with an orbit on which G acts (faithfully). The
# following method does just that. The caller is responsible for supplying
# a suitable vector v (the shorter the orbit, the faster the computations).
#
SizeViaOrbit := function (G, v)
	local orbit, phi, size;

	orbit := Orbit(G, v);;
	orbit := ShallowCopy(orbit);; Sort(orbit);
	phi := ActionHomomorphism(G, orbit);
	size := Size(Image(phi));
	Unbind(phi); Unbind(orbit);
	return size;
end;;

#
# Find a nice permutation presentation of a matrix group by its
# action on the orbit of a given vector.
#
IsomorphismPermGroupViaOrbit := function (G, v)
	local orbit, phi;

	Print("  Computing orbit...\n");
	orbit := Orbit( G, v );;
	orbit := ShallowCopy(orbit);;
	Sort(orbit);
	#
	Print("  Computing permutation group ...\n");
	phi := ActionHomomorphism( G, orbit, "surjective" );;
	#
	Print("  Verifying that perm group is isomorphic image of the stabilizer...\n");
	Assert(0, IsSurjective(phi) and Size(Image(phi)) = Size(G));
	#
	return phi;
end;;

#
# Find relators of a finite presentation on FG_gens of the group H generated
# by H_gens, where each element of H_gens is mapped to the corresponding
# element of FG_gens.
#
GroupRelatorsViaOrbit := function( H, H_gens, FG_gens )
	local psi, Fgens, Frels;

	Print("  Determining finite presentation...\n");
	psi := IsomorphismFpGroupByGenerators( H, H_gens );;

	# Determine the relators
	Fgens := FreeGeneratorsOfFpGroup(Image(psi));
	Frels := RelatorsOfFpGroup(Image(psi));

	# Verify all relations actually hold in H
	Assert(0, ForAll( Frels, r->MappedWord( r, Fgens, H_gens ) = One(H) ));
	return List( Frels, r->MappedWord( r, Fgens, FG_gens ) );
end;;

#
# Auxillary function used by G2 and Sp6 below. Computes a list of
# generators for G2 in the representation we need.
#
G2_Sp_helper := function (q, TMP)
	local F, gens, g, g2, x, i, j, G;

	Assert(0, IsEvenInt(q));

	F := GF(q);
	gens := [];

	for g in GeneratorsOfGroup(TMP(3,F)) do
		g2 := TransposedMat(g)^-1;
		x := NullMat(6,6,F);
		for i in [1..3] do
			for j in [1..3] do
				x[i][j] := g[i][j];
				x[3+i][3+j] := g2[i][j];
			od;
		od;
		Add(gens, x);
	od;

	x := [
		[ 1,  0,  0,  1,  0,  0],
		[ 0,  1,  0,  0,  0,  0],
		[ 0,  0,  1,  0,  0,  0],
		[ 0,  0,  0,  1,  0,  0],
		[ 0,  0,  1,  0,  1,  0],
		[ 0, -1,  0,  0,  0,  1] ] * One(F);
	Add(gens, x);

	return Group(gens);
end;;

#
# Return a suitable representation of the group G2(q),
# where q must be a power of two.
#
G2 := function(q)
	local G;
	G := G2_Sp_helper(q, SL);
	#SetSize(G, q^6*(q^6-1)*(q^2-1));
	return G;
end;;

#
# Return a suitable representation of the group Sp(6,q),
# where q must be a power of two.
#
Sp6 := function(q)
	local G, x;
	G := G2_Sp_helper(q, GL);
	if q = 2 then
		x := [
			[ 1,  0,  0,  1,  0,  0],
			[ 0,  1,  0,  0,  1,  0],
			[ 0,  0,  1,  0,  0,  1],
			[ 0,  0,  0,  1,  0,  0],
			[ 0,  0,  0,  0,  1,  0],
			[ 0,  0,  0,  0,  0,  1] ] * One(GF(q));
		G := ClosureGroup(G, x);
	fi;
	#SetSize(G, q^6*(q^6-1)*(q^2-1));
	return G;
end;;
\end{lstlisting}

}{} 

\noindent Authors' Addresses: 
 
{\footnotesize 
\medskip 
\noindent \hspace{-2cm}\begin{tabular}{lll}Ralf Gramlich, Max Horn & Antonio Pasini & Hendrik Van Maldeghem \\ 
TU Darmstadt & Universit\`a di Siena & Ghent University\\ 
FB Mathematik / AG 5 & Scienze Matematiche e Informatiche & Pure Mathematics and Computer Algebra\\ 
Schlo\ss gartenstra\ss e 7 & Pian dei Mantellini 44 & Krijgslaan 281, S22\\ 
64289 Darmstadt & 53100 Siena & 9000 Gent\\ 
Germany & Italy & Belgium \\ 
{\tt gramlich@mathematik.tu-darmstadt.de}\phantom{lala}& {\tt pasini@unisi.it} & {\tt hvm@cage.ugent.be} \\ 
{\tt mhorn@mathematik.tu-darmstadt.de} \\
\\
First author's alternative address: \\
Ralf Gramlich \\
The University of Birmingham \\
School of Mathematics \\
Birmingham \\
B15 2TT \\
United Kingdom \\
{\tt ralfg@maths.bham.ac.uk}
\end{tabular} 
}

\end{document}